\documentclass[11pt]{amsart}

\usepackage{amssymb,upgreek}
\usepackage{subcaption}
\usepackage{url}
\usepackage{bm}
\usepackage[left=1in,top=1in,right=1in,bottom=1in,head=.2in]{geometry}
\usepackage[titletoc]{appendix}
\usepackage[percent]{overpic}

\usepackage[textsize=scriptsize]{todonotes}
\usepackage{color}
\usepackage{mathtools} 

\usepackage{fancyhdr}
\usepackage{mathrsfs}
\usepackage[cmtip,all]{xy}

\usepackage{hyperref}
\usepackage{enumerate}

\swapnumbers
\newtheorem{theorem}{Theorem}[section] 
\newtheorem{lemma}[theorem]{Lemma}

\newtheorem{corollary}[theorem]{Corollary}
\newtheorem*{theorem*}{Theorem}

\newtheorem*{fcthm*}{Finite Cork Theorem}
\newtheorem*{ccthm*}{Cork Consolidation Theorem}

\newtheorem*{thm2*}{Light Bulb Theorem for disks (version ii)}
\newtheorem*{lbthm*}{Generalized 4D Lightbulb Theorem}
\newtheorem*{2lbthm*}{Generalized 4D Lightbulb Theorem (restated)}
\newtheorem*{icthms*}{Infinite Cork Theorems}
\newtheorem*{aclemma*}{\ac-Lemma}
\newtheorem*{mclemma*}{Multicork Lemma}
\newtheorem*{multicorktheorem*}{Multicork Theorem}
\newtheorem*{lemma*}{Lemma}
\newtheorem*{corollary*}{Corollary}

\newcommand{\thistheoremname}{}
\newtheorem{genericthm}[theorem]{\thistheoremname}

\theoremstyle{definition}
\newtheorem{definition}[theorem]{Definition}
\newtheorem{question}[theorem]{Question}
\newtheorem{remark}[theorem]{Remark}

\newtheorem{example}[theorem]{Example}

\newtheorem*{remark*}{Remark}
\newtheorem*{definition*}{Definition}
\newtheorem*{remarks*}{Remarks}
\newtheorem*{addenda*}{Addenda}

\newcommand{\pf}{\vskip-5pt \vskip-5pt \proof}

\newcommand{\fig}[3]{\begin{figure}\includegraphics[height=#1pt]{#2}#3\end{figure}}

\newcommand{\bit}[1]{\textbf{\textit{#1}}} 

\newcommand{\pt}{\textup{pt}}

\newcommand{\sto}{\!\!\xymatrix@C=1em{{}\ar@{~>}[r]&{}}\!\!}

\newcommand{\interior}{\textup{int}}
\newcommand{\ac}{\textup{AC}}

\newcommand{\fq}{\operatorname{fq}}
\newcommand{\dax}{\operatorname{dax}}
\newcommand{\Dax}{\operatorname{Dax}}

\newcommand{\items}{\begin{itemize}[leftmargin=25pt,rightmargin=5pt]
  \setlength\itemsep{2pt}}
\newcommand{\stopitems}{\end{itemize}}

\address{North Carolina School of Science and Mathematics, 901 Burkemont Avenue,
Morganton, NC 28655}
\email{hannah.schwartz@ncssm.edu} 

\begin{document}

\title{A 4-dimensional light bulb theorem for disks}
\author{Hannah R. Schwartz}

\begin{abstract}
We give a $4$-dimensional light bulb theorem for properly embedded disks, generalizing recent work of Gabai and Kosanovi\'c-Teichner in certain contexts, and extending the $4$-dimensional light bulb theorem for $2$-spheres due to Gabai and Schneiderman-Teichner. In particular, we provide conditions under which homotopic disks properly embedded in a compact $4$-manifold $X$ with a common dual in the interior of $X$ are smoothly isotopic rel boundary. We also provide a new geometric interpretation of the Dax invariant, to aid in its computation. 
\end{abstract}

\maketitle

\vskip-.4in
\vskip-.4in

\parskip 2pt

\setcounter{section}{-1}

\parskip 2pt

\section{Introduction and Motivation}

We are interested in an isotopy classification of pairs of homotopic disks $D_0$ and $D_1$ properly embedded with a common dual in a smooth, compact $4$-manifold $X$ with non-empty boundary $\partial X$. As has become rather standard in the literature, ``common dual" refers to a $2$-sphere embedded in $X$ with trivial normal bundle that intersects both $D_0$ and $D_1$ transversally in a single point. Recent work of Gabai \cite{dave:LBT} and Schneiderman-Teichner \cite{st} completely characterizes the conditions under which homotopic spheres $S, T \subset X$ with a common dual are isotopic. Namely, there is a smooth isotopy between $S$ and $T$ whenever the concordance invariant defined in Theorem $10.5A$ of Freedman and Quinn \cite{freedman-quinn:4-manifolds} vanishes. Indeed, the author in \cite{me} constructed infinitely many explicit examples of pairs of spheres for which this invariant is non-zero. 

Recent work of Gabai \cite{dave:LBL} and Kosanovi\'c-Teichner \cite{KT} gives analogous ``light bulb theorems" for homotopic properly embedded disks with a common dual in the \emph{boundary} $\partial X$ and vanishing Dax invariant, an isotopy invariant for pairs of homotopic properly embedded disks based on homotopy theoretic work of Dax \cite{dax}. Gabai \cite{dave:LBL} was first to realize the Dax invariant as an isotopy obstruction for disks; later Kosanovi\'c and Teichner \cite{KT} proved that the Dax invariant completely classifies smooth isotopy of disks in this setting. 
 
Our main results furthers the isotopy classification of homotopic, properly embedded disks with a common dual in the \emph{interior} of $X$. We first answer Question $5.7$ posed by Gabai in \cite{dave:LBL}. 

\begin{lemma}\label{Dgiso}
Let $D$ be a properly embedded disk in $X$, with dual $G$ in the interior of $X$. For any $g \in \pi_1(X)$ there is a smooth isotopy rel boundary between the self-referential disks $D_g$ and $D_{{\bar{g}}}$ supported away from their common dual $G$. 
\end{lemma}

Background on self-referential disks is given in  Definition \ref{selfref} and Figure \ref{disks}. This lemma is critical to the proof of Theorem \ref{main} below, our version of a 4D light bulb theorem for disks.

\begin{theorem} \label{main} 
Let $D_0$ and $D_1$ be properly embedded disks in a $4$-manifold $X$. Suppose that $D_0$ and $D_1$ are homotopic rel boundary, with vanishing Freedman-Quinn invariant $\rm{fq}(D_0, D_1)$. If $D_0$ and $D_1$ have a common $\pi_1$-negligible dual sphere $G$, then there is a smooth isotopy from the disk $D_1$ to a self-referential disk $(D_0)_\omega$ for some $\omega \in \mathbb{Z}[\pi_1(X)]$ with $\mathcal D(\omega) = \rm{Dax}(D_0, D_1)$ in the quotient $\mathbb{Z}[\pi_1(X) \setminus 1]/ \rm{dax}(\pi_3(X))$. This isotopy can be assumed to be rel boundary and supported away from the dual $G$. 
\end{theorem} 

In other words, homotopic disks with a common dual differ up to isotopy only by a collection of self-referential tubes that depend on the Dax invariant of the pair of disks. In particular, if the manifold $X$ is simply-connected, then homotopy implies isotopy just as in other settings.
 
This main theorem generalizes Theorem $2.5$ of Gabai \cite{dave:LBL}, in which the dual $G$ is required to be in the boundary $\partial X$. In fact, Gabai states the need for our $\pi_1$-negligibility assumption on the dual sphere in this more general case, in his Remark $2.7$(iii). A similar condition on the dual sphere is also required in Gabai's light bulb theorem for higher genus surfaces in Theorem $9.7$ of \cite{dave:LBT}. The question of homotopy versus isotopy for properly embedded disks without this $\pi_1$ condition on the dual is still open; see Question \ref{q}.

The final section \ref{appendix} includes a new geometric formulation of the Dax invariant, along with a sample computation. 

\smallskip
\noindent
\bit{Acknowledgements.}
The author would like to thank Dave Gabai and Peter Teichner for their guidance, encouragement, and interest in these results, as well as Rob Schneiderman for helpful discussions and thoughtful comments from the referees. A large part of this project was completed while the author was visiting the Max Planck Institute for Mathematics; she thanks the institute for providing a welcoming and collaborative research environment, as well as Voigt for providing the baked goods necessary to maintain her mental aptitude. The rest was completed at Princeton University and the North Carolina School of Science and Mathematics, where the author is supported by NSF grant DMS-1502525.

\section{Preliminaries}

In this section, we introduce some background necessary for Sections \ref{2} and \ref{3}. Although the disks in our final setting will always have a dual sphere, we begin below in a more general context without a dual until it is introduced formally in Section \ref{dualsec}. We work throughout in the smooth and oriented category, with all homotopies and isotopies between disks taken rel boundary, and all maps (including homotopies) generic unless otherwise stated. In general, our convention will be to refer to double points of immersions using the letters $p$ and $q$, while their pre-images will be denoted by the letters $x$ and $y$.

\subsection{Regular homotopies} \label{reghom} 
Let $D_0$ and $D_1$ be disks properly embedded in a compact $4$-manifold $X$ with boundary, related by a homotopy $h:I \times D^2 \to X^4$. 

\fig{100}{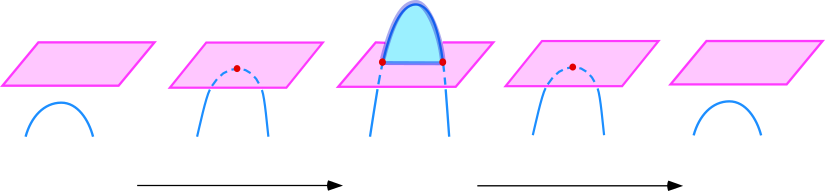}{
\put(-221,78){$W$}
\put(-221,59){\small $\omega$}
\put(-204,93){\small $\omega^*$}
\put(-196,63){\small $-$}
\put(-243,63){\small $+$}
\put(-168,-10){\small{Whitney move}}
\put(-340,-10){\small{Finger move}}
\caption{$3$-dimensional ``slices" of the local model of the disk throughout the regular homotopies from Definition \ref{fwmovedef}.}
\label{fw}}

\begin{definition}\label{fwmovedef} 
By Smale \cite[Theorem D]{smale}, we can and will always assume that $h$ is a \bit{regular homotopy}, i.e.~a homotopy through smooth immersions. There are only finitely many times during such a homotopy at which the immersed disk is not self-transverse: when pairs of oppositely signed double points in its interior are either introduced or cancelled. The local model for the regular homotopy shown in Figure \ref{fw} that introduces pairs of oppositely signed double points is called a \bit{finger move}, and its inverse homotopy is called a \bit{Whitney move}. A Whitney move is supported in the regular neighborhood of an embedded \bit{Whitney disk} $W$ with boundary the union of two \bit{Whitney arcs} $\omega$ and $\omega^*$ embedded in the immersed disk that pair the double points. Figure \ref{fw} illustrates both a finger and Whitney move; note that these diagrams are truly \emph{local models} of the homotopy rather than simply schematics, as are all of the illustrations of regular homotopies that follow. To explicitly define a Whitney move in local coordinates requires that the normal disk bundle of the Whitney disk be framed ``compatibly" with respect to its boundary on the immersed disk. See \cite{freedman-quinn:4-manifolds} as well as Casson's lectures in \cite{casson} for more exposition. 
\end{definition} 

By general position, the arc along which any finger move is supported can be isotoped away from all Whitney disks. It follows that there is an ambient isotopy of $X$ rel boundary taking the homotopy $h$ to one that factor into a collection of simultaneous finger moves taking the disk $D_0$ to an immersed disk $h_{1/2}(D^2):= \Sigma$, followed by a collection of simultaneous Whitney moves taking $\Sigma$ to the embedded disk $D_1$. There are therefore \emph{two} distinguished sets of Whitney disks $\mathcal W$ and $\mathcal V$ pairing the double points of $\Sigma$, along which Whitney moves give the embeddings $D_0$ and $D_1$, respectively. 

\begin{definition} \label{wdisks}
We refer to the immersed disk $\Sigma$ above as the \bit{middle level} of the regular homotopy, and the Whitney disks in the sets $\mathcal W$ and $\mathcal V$ as the \bit{descending and ascending Whitney disks}, with boundaries the \bit{descending and ascending Whitney arcs}. 
\end{definition}

\subsection{Dual spheres and tubed disks} \label{dualsec}

\begin{definition} \label{dual}
A smoothly embedded $2$-sphere $G \subset X$ is said to be \bit{dual} to an embedded surface $F \subset X$ if
\begin{enumerate}
\item $G$ has trivial normal bundle in $X$, and
\item $G$ intersects $F$ transversally in a single point. 
\end{enumerate} 
The dual $G$ is called \bit{$\pi_1$-negligible} if the map $\pi_1(X-G) \to \pi_1(X)$ induced by inclusion is an isomorphism, or equivalently, if $G$ has an immersed dual sphere in $X$.  
\end{definition}

The existence of dual spheres was most notably used by Gabai \cite{dave:LBT} and Schneiderman-Teichner \cite{st} to convert regular homotopies into smooth isotopies. Most recently in \cite{kst}, the notion of a dual has been useful even up a dimension, to produce isotopies of 2-spheres admitting algebraic dual 3-spheres in a 5-manifold. 

It is helpful in many of our arguments to mimic portions of both \cite{dave:LBT} and \cite{st}, especially by working in the context of ``tubed surfaces" introduced by Gabai in \cite{dave:LBT}. Many of Gabai's constructions were originally in the context of either spheres with duals, or properly embedded disks with duals in the boundary. Nonetheless, we apply many of Gabai's constructions here, careful to note any instances where this change of context is relevant. Remark $2.7$ of Gabai \cite{dave:LBL} is enlightening in this vein, as he highlights where and how such assumptions are needed in his previous work. 

\fig{300}{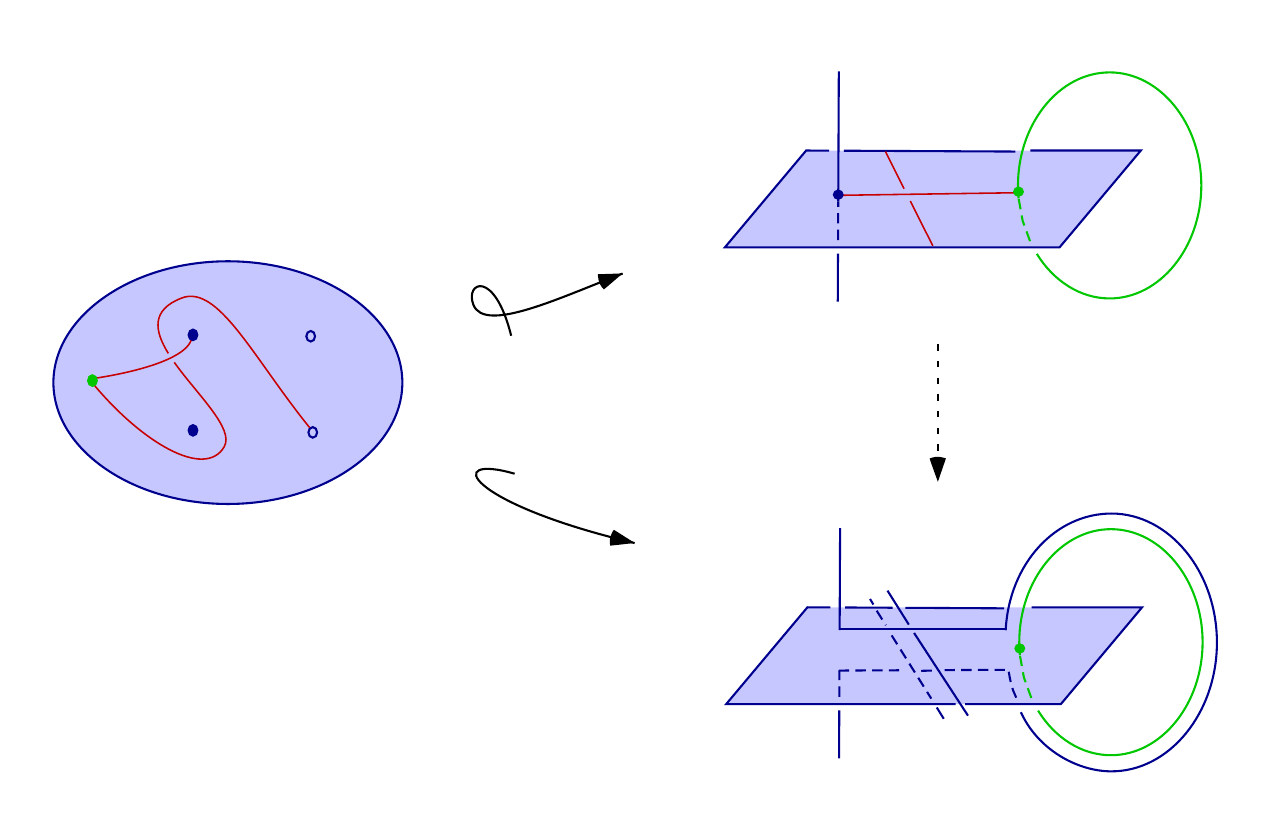}{
\put(-387,98){$D^2$}
\put(-120,233){$\alpha_1$}
\put(-150,250){$\alpha_2$}
\put(-185,196){$\Sigma$}
\put(-195,29){$(\Sigma,\alpha)$}
\put(-52,37){\small{$G$}}
\put(-52,203){\small{$G$}}
\put(-110,150){\small{Norman trick}}
\caption{An immersed disk $\Sigma \subset X$ with two double points and dual $G$ (top), together with a collection $\alpha=\{\alpha_1, \alpha_2\}$ of tubing arcs in $\Sigma$ connecting each double point to the point $\Sigma \pitchfork G$. As in Definition \ref{tubeddisk}, this determines an embedded tubed disk $(\Sigma,\alpha)$ (bottom). It is often convenient to picture isotopies between tubed disks as a motion of the pre-images of the tubing arcs in $D^2$ (left).}
\label{tubes}}

\begin{definition} \label{tubeddisk} 
An  immersed disk $\Sigma \subset X$ with $n$ double points and dual $G$, together with a collection $\alpha=\{\alpha_1, \dots, \alpha_n\}$ of (possibly immersed) \bit{tubing arcs} in $\Sigma$ connecting each double point to the point $\Sigma \pitchfork G$, determines an embedded \bit{tubed disk $(\Sigma, \alpha)$} which also has dual $G$, constructed by using the Norman trick \cite{norman} to tube each double point of $\Sigma$ over $G$ along the tubing arc $\alpha$ as in Figure \ref{tubes}. If the tubing arcs in $\alpha$ are immersed, then the isotopy class of the resulting tubed disk $(\Sigma,\alpha)$ is not unique until an ordering of the arcs (i.e. an overstrand and an understrand) is chosen at each crossing to determine the relative radius of the tubes running along those arcs in the normal disk bundle of $\Sigma$, as shown in Figure \ref{tubes}. A diagram on the disk $D^2$ of the pre-image of the tubing arcs, together with their endpoints and crossing information, is called a \bit{tube diagram} for the tubed disk. 
\end{definition}

Isotopies of tubed disks which only move the ``tubes" are convenient to visualize, especially from the point of view of the pre-image: the pre-images of the tubing arcs in $D^2$ determine how the tubes run along the embedded tubed disk in $X$. In particular, Gabai's isotopies of tubed surfaces from \cite[Figure 5.6]{dave:LBT} will come in handy.

\fig{280}{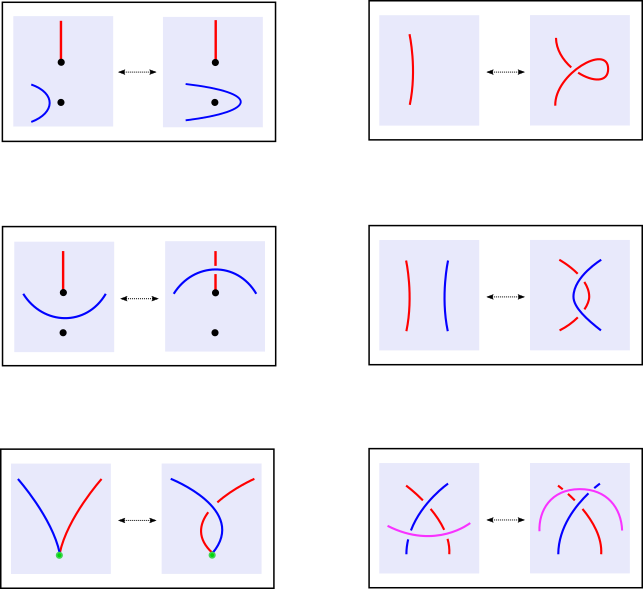}{
\put(-265,200){\small{\bf Isotopy A}}
\put(-265,95){\small{\bf Isotopy B}}
\put(-71,200){\small{\bf RI}}
\put(-74,95){\small{\bf RII}}
\put(-74,-12){\small{\bf RIII}}
\put(-270,-12){\small{\bf Re-ordering}}
\caption{Isotopies of a tubed disk which move only the tubes and leave the rest of the disk fixed, as in Definition \ref{AB}. These isotopies, drawn from the perspective of the tube diagram, are identical to those defined in \cite[Figure 5.6]{dave:LBT}. Re-ordering is also described by Schneiderman-Teichner in Figures $14$ and $15$ of \cite{st}}
\label{moves2}}

\fig{140}{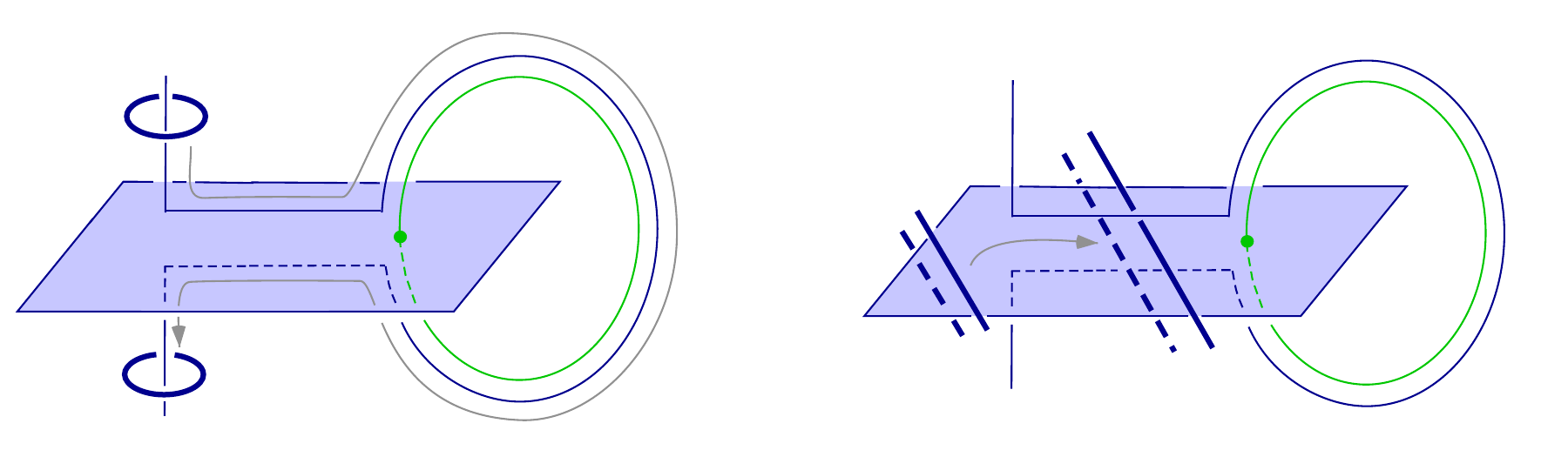}{
\put(-410,-10){\small{\bf Isotopy A}}
\put(-150,-10){\small{\bf Isotopy B}}
\put(-80,97){\small{$G$}}
\put(-340,97){\small{$G$}}
\put(-210,77){\small{$\tau$}}
\put(-159,101){\small{$\tau'$}}
\put(-452,101){\small{$\tau$}}
\put(-454,23){\small{$\tau'$}}
\caption{Three dimensional cross sections of the disks during the tube isotopies A and B from Definition \ref{AB} and Figure \ref{moves2}: Isotopy A swings the tube $\tau$ across the dual $G$ to a new tube $\tau'$ as indicated by the arrows, and isotopy B moves the tube $\tau$ to a tube $\tau'$ that is nested ``inside" of another tube.}
\label{moves}}

\begin{definition} \label{AB} Any composition of the isotopies of tubed disks illustrated in Figure \ref{moves2}, as well as \cite[Figure 5.6]{dave:LBT}, will be referred to as a \bit{tube isotopy}. These isotopies are ``local" in the sense that they move only one tube at a time and leave the rest of the disk fixed. The support of each tube isotopy is also completely contained in a tubular neighborhood of the original surface and its dual sphere. The isotopies in Figure \ref{moves2} are drawn from the perspective of tube diagrams, whereas Figure \ref{moves} shows the $3$-dimensional cross-sections of the tubed disk during isotopies $A$ and $B$.  
\end{definition} 

Another benefit of dealing with tubed surfaces is the ability to ``shadow" a homotopy, in the following sense due to Gabai \cite[Section 5]{dave:LBT}. Consider a regular homotopy $h$ between embedded disks $D_0$ to $D_1$ with a common dual $G$ and middle level $\Sigma$. Suppose that $h$ is supported away from $G$. In this case, the sphere $G$ is dual to $\Sigma$, and disjoint from both the descending Whitney disks $\mathcal W=\{W_1, \dots, W_n\}$ and ascending Whitney disks $\mathcal V=\{V_1, \dots, V_n\}$ from Definition \ref{wdisks}. 

\begin{definition} \label{shadow}
For $i=1, 2, \dots, n$, let $\omega_i \subset \partial W_i$ denote a choice of one of the two descending Whitney arcs on the boundary of $W_i$, with endpoints $p_i^\pm$ the $i^{th}$ pair of positive and negative double points of $\Sigma$. First, choose (possibly immersed) arcs $\nu_i$ connecting the interior of $\omega_i$ to the point $\Sigma \pitchfork G$ in the complement of all Whitney arcs of $\mathcal W$. Then, define tubing arcs $\alpha= \{\alpha^\pm_1, \dots, \alpha^\pm_n\}$ such that $\alpha_i^\pm$ runs first from $p_i^\pm$ along $\omega_i$ to $\nu_i$, and then along $\nu_i$ to $\Sigma \pitchfork G$. The tubed disk $(\Sigma,\alpha)$, illustrated in Figure \ref{shadowpic}, is called a \bit{descending shadow of $h$}. A tubed disk where instead \emph{ascending} Whitney arcs are used to construct the tubing arcs is called an \bit{ascending shadow of $h$}. 
\end{definition}

By Gabai \cite[Lemma 5.1]{dave:LBT}, the disk $D_0$ is smoothly isotopic (rel boundary and rel the dual $G$) to the descending shadow $(\Sigma, \alpha)$ from Definition \ref{shadow}, while the disk $D_1$ is smoothly isotopic to the ascending shadow of $h$. The term ``shadow" references the fact that an explicit isotopy from $D_0$ (respectively $D_1$) to the ascending shadow (respectively, descending shadow) can be defined by following each finger move to $\Sigma$ except for swinging around the dual $G$ to avoid creating double points. Note that if the descending and ascending shadows are tube isotopic, then the disks $D_0$ and $D_1$ are isotopic rel boundary as well. 

\fig{320}{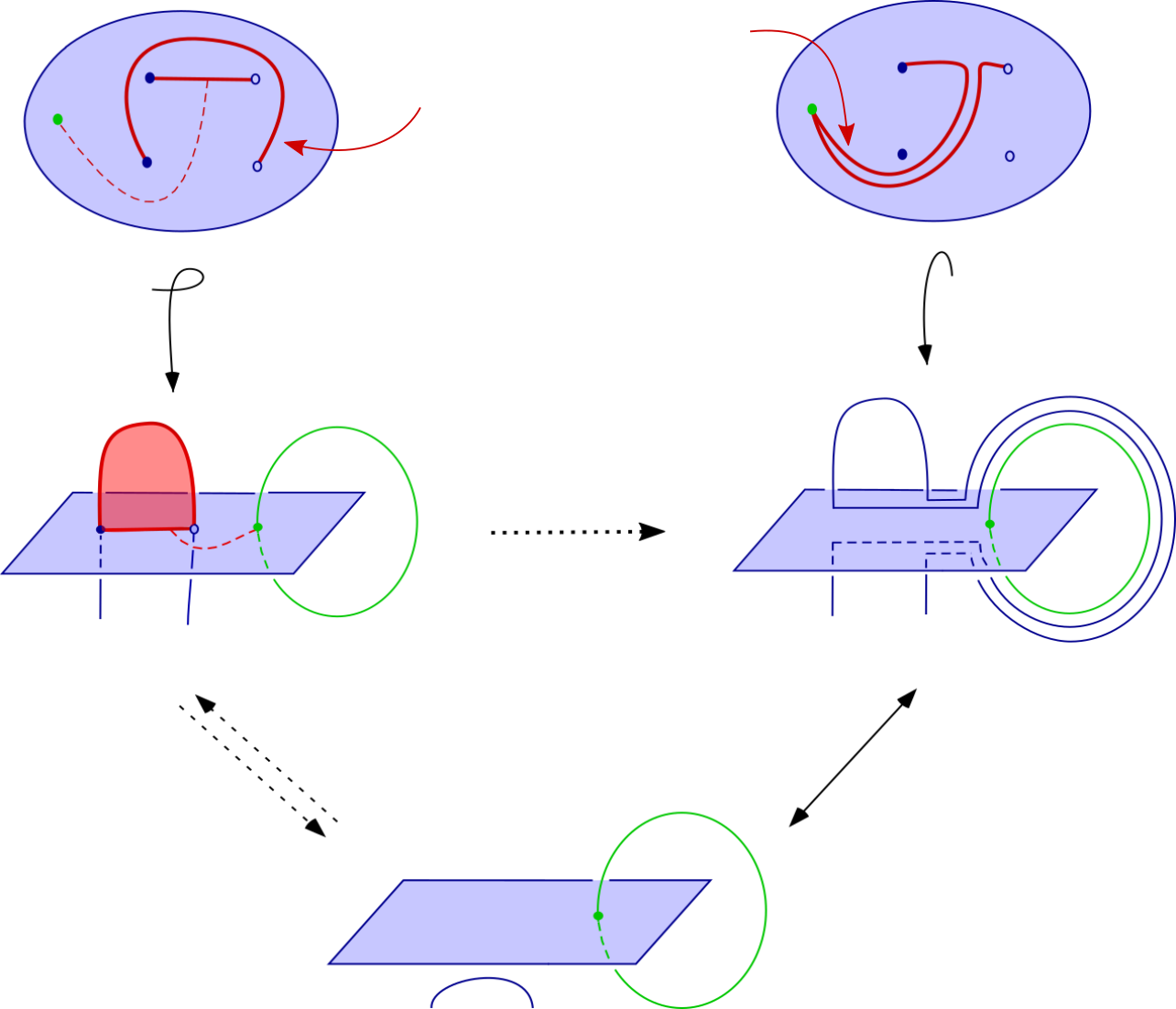}{
\put(-254,300){\small{Pre-image of}}
\put(-255,290){\small{Whitney arcs}}
\put(-193,318){\small{Pre-image of}}
\put(-192,308){\small{tubing arc $\alpha$}}
\put(-144,113){\small{$(\Sigma, \alpha)$}}
\put(-270,128){\small{$G$}}
\put(-333,168){\small{$W_i$}}
\put(-365,126){\small{$\Sigma$}}
\put(-220,160){\small{Norman trick}}
\put(-286,84){\small{f.m.}}
\put(-317,70){\small{w.m.}}
\put(-93,75){\small{isotopy}}
\put(-153,126){\small{``shadow"}}
\put(-36,130){\small{$G$}}
\put(-159,4){\small{$G$}}
\put(-260,2){\small{$D_0$}}
\caption{A descending shadow $(\Sigma, \alpha)$ of the homotopy $h$ in Definition \ref{shadow}. By Gabai \cite[Lemma 5.1]{dave:LBT}, the tubed disk $(\Sigma, \alpha)$ is isotopic to $D_0$ (as indicated on the bottom right). The top row shows the disks in the domains of the immersions, decorated by the pre-images of the double points of $\Sigma$, as well as pre-images of the tubing and Whitney arcs.}
\label{shadowpic}}

\subsection{Self-referential disks and the Dax invariant}  \label{daxsection}

Let $D$ be a properly embedded disk in a 4-manifold $X$.  For fundamental group calculations, assume that the basepoint $b$ of $\pi_1(X)$ is chosen in $\partial D \subset \partial X$. 

\begin{definition}\label{selfref} 
Let $B_1, \dots, B_n$ be $3$-balls disjointly embedded in $X-D$, and $\eta_1, \dots, \eta_n$ disjointly embedded arcs with interiors in $(X - D) - \bigcup_i \partial B_i$ such that 
\begin{enumerate}
\item each $\eta_i$ connects a point $z_i \in \interior(D)$ to a point $\widehat z_i \in \partial B_i$, and
\item $|\eta_i \pitchfork \interior(B_j)| = \delta_{ij}$.
\end{enumerate}

Tubing the disk $D$ to each sphere $\partial B_i$ along the arc $\eta_i$ gives the \bit{self-referential} embedded disk $D_\omega \subset X$ for the word $$\omega= \sum_{i=1}^n\epsilon_i g_i \in \mathbb{Z}[\pi_1(X)]$$ with $\epsilon_i = \pm 1$ the sign of the intersection point $z_i=\eta_i \pitchfork \interior(B_i)$, and $g_i \in \pi_1(X)$ the homotopy class of a based loop which travels (in order) along: {\bf (i)} a path $\nu \subset D$ from the basepoint to $z_i$, {\bf (ii)} the path $\eta_i$, {\bf (iii)} a path in $B_i$ from $\widehat z_i$ to $z_i$, {\bf (iv)} the second ``half" of the path $\eta_i^{-1}$ from $\widehat z_i$ to $z_i$, and finally {\bf (v)} the path $\nu^{-1}$. See Figure \ref{disks}. 
\end{definition}

The construction of self-referential disks (and ``self-referential homotopies", defined below) is due to Gabai \cite{dave:LBL} and based on the \emph{spinning} construction of Budney \cite{budney}. A similar construction can be found in \cite[Section $5.2$]{KT}, see Figure $18$ in particular. 

\begin{remark}
If the disk $D$ has a dual sphere $G$, as in Definition \ref{dual}, we will always assume that the balls $B_i$ and arcs $\eta_i$ are chosen away from $G$. Thus the self-referential disk $D_\omega$ has the same dual $G$. 
\end{remark}

Lemma $2.15$ of Gabai \cite{dave:LBL} shows that if $D$ has a dual, then the isotopy class of a self-referential disk $D_\omega$ depends only on the word $\omega$ and not on individual cancelling terms $+g_i$ and $-g_i$ appearing in the sum. Although Gabai is working with a dual sphere in the boundary of the ambient $4$-manifold, his proof of this lemma does not depend on this fact (this can be seen from Figure $7$ of \cite{dave:LBL}). Indeed, the argument works with any embedded dual in the interior of the ambient manifold, needed so that the "lightbulb isotopy" as in \cite{dave:LBT} can be performed. 

\begin{definition} \label{selfrefhom}
Given a self-referential disk $D_\omega$ for $\omega \in \mathbb{Z}[\pi_1(X)]$, there is a regular \bit{self-referential homotopy $h_\omega$} from $D_\omega$ to the disk $D$ that ``unlinks"  each sphere $\partial B_i$ from the tube running along the arc $\eta_i$ from Definition \ref{selfref}. In particular, one finger move and one Whitney move can be used to unlink each tube, as shown in Figure \ref{Dg} when $\omega=g$ and $\omega=\bar g$ for any $g \in \pi_1(X)$; these cases will be of particular significance to us. In general, a self-referential homotopy $h_\omega$ consists of performing a sequence of regular homotopies from Figure \ref{Dg} on each tube of $D_\omega$.
\end{definition}

\fig{170}{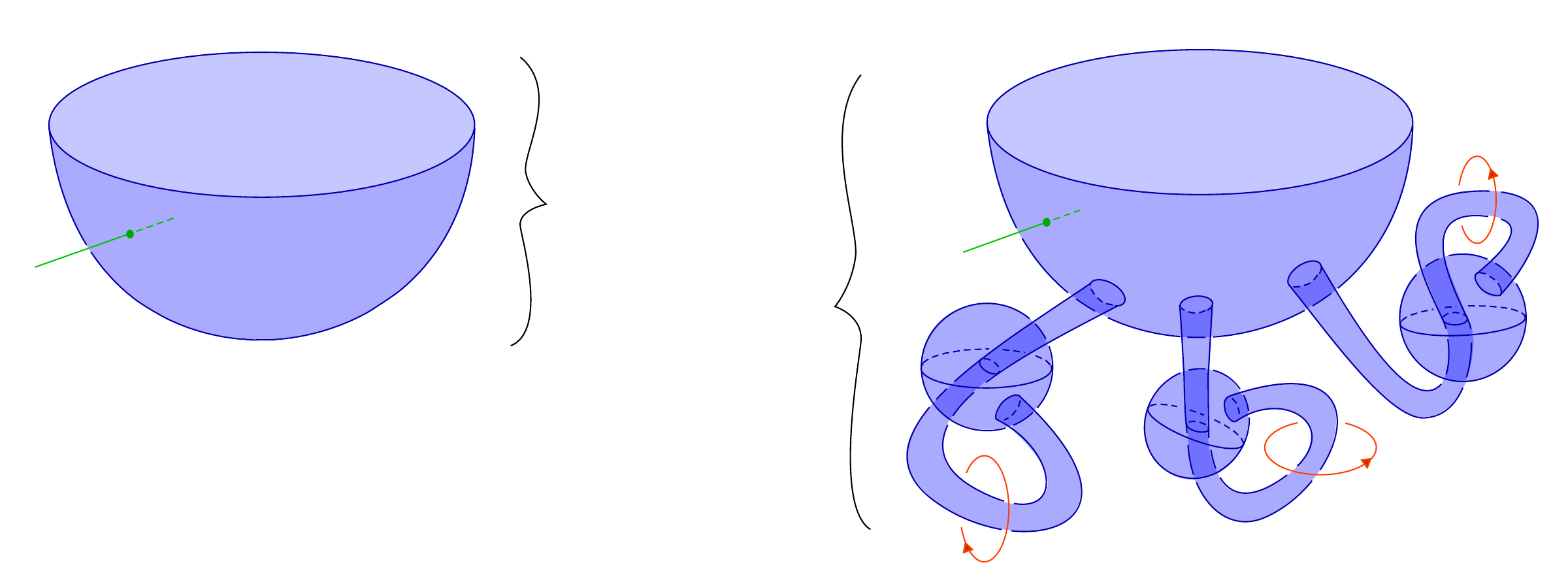}{
\put(-455,88){\small{$G$}}
\put(-188,93){\small{$G$}}
\put(-290,107){\small{$D$}}
\put(-231,79){\small{$D_\omega$}}
\put(-172,0){\small{$g_1$}}
\put(-52,38){\small{$g_2$}}
\put(-29,129){\small{$g_3$}}
\caption{An example of a self-referential disk $D_\omega$ as in Definition \ref{selfref}, with word $\omega \in \mathbb{Z}[\pi_1(X)]$ equal to the sum $g_1 + g_2 - g_3$ for $g_i \in \pi_1(X)$. In our schematic, loops representing the group element $g_i$ are portrayed as linking the corresponding oriented red circle $+1$ times.}
\label{disks}}

\fig{500}{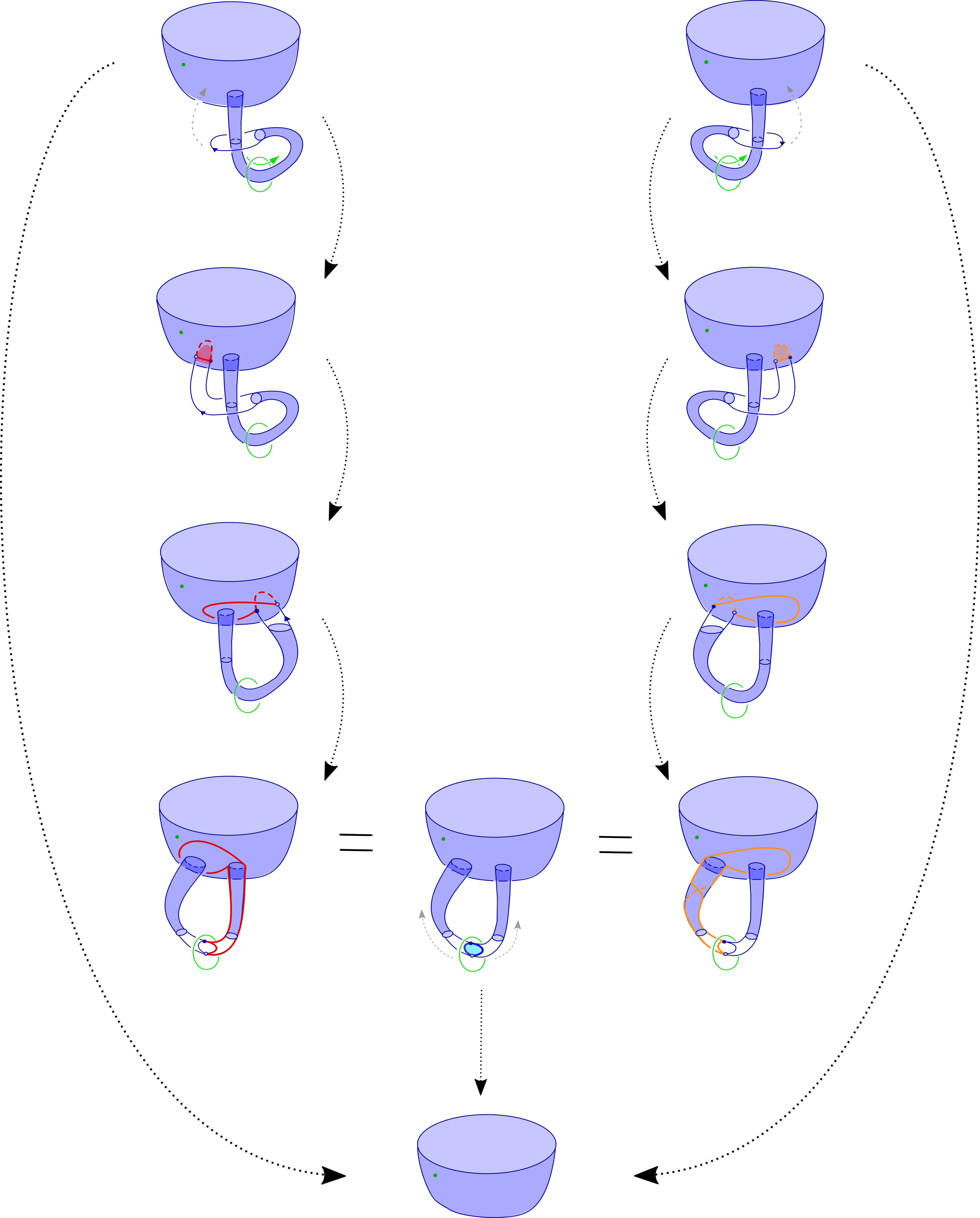}{
\put(-270,475){\small{$D_g$}}
\put(-140,475){\small{$D_{\bar g}$}}
\put(-68,130){\small{$\Sigma$}}
\put(-173,130){\small{$\Sigma$}}
\put(-288,130){\small{$\Sigma$}}
\put(-155,425){\small{f.m.}}
\put(-255,425){\small{f.m.}}
\put(-420,260){\small{$h_g$}}
\put(10,260){\small{$h_{\bar g}$}}
\put(-177,-5){\small{$D$}}
\put(-200,74){\small{w.m.}}
\caption{Local models of self-referential homotopies $h_g$ and $h_{\bar g}$, as in Definition \ref{selfrefhom}, from the self-referential disks $D_g$ (upper left) and $D_{\bar g}$ (upper right) to the disk $D$ (bottom center). The green circle in each diagram indicates the group element $g \in \pi_1(X)$, where a loop representing $g$ links the circle in the direction shown in the top two diagrams. Both $h_g$ and $h_{\bar g}$ consist of a finger move supported away from the dual $G$ to the same middle level $\Sigma$, followed by a Whitney move along the Whitney disk shaded blue (fourth row center) taking $\Sigma$ to the disk $D$. The initial finger moves of $D_g$ and $D_{\bar g}$ are inverse to Whitney moves along the Whitney disks shaded in the second row. The Whitney arcs on the boundary of these Whitney disks are drawn in red and orange, respectively, on the middle level $\Sigma$. The ``twist" in the Whitney arc of for $h_{\bar g}$ relative to the Whitney arc for $h_g$ is needed so that the positive and negative double point align in the immersion $\Sigma$. All unlabelled arrows in the diagram denote motions of $\Sigma$ induced by ambient isotopies of $X$. }
\label{Dg}}

\begin{definition}\label{dax}
Let $h: I \times D^2 \to X^4 $ be a homotopy rel boundary between properly embedded disks. From now on, we shall fix a parametrization of the disk as a product $D^2 = J \times K$ of two (oriented) unit intervals. Then, for each ``time" $t \in I$ during the homotopy, the immersion of the disk $h_t: J \times K \to X^4$ can be thought of as a homotopy of the interval $K$ in $X$ fixing its boundary endpoints in $\partial X$. By work of Dax \cite{dax}, the smooth map $h$ may be perturbed (rel boundary) so that the map $\widehat h: I \times J \times K \to I \times J \times X$ sending $(t, s, r) \mapsto (t, s, h_t(s, r))$ is an immersion that has no triple points, and only finitely many transverse double points $p_1, \dots, p_n$ with pre-images $x_i, y_i \in I \times J \times K$ such that $x_i$ has $K$-coordinate strictly less that that of $y_i$. 

A unique group element $g_i \in \pi_1(X)$, represented by the based loop $\gamma_i$ in Figure \ref{daxpic}, can also be assigned to each double point. For simplicity, assume that the double points are labelled so that $g_i$ is non-trivial for $i=1,  \dots, m$ and trivial otherwise. A well-defined sign $\epsilon_i = \pm 1$ for each double point $p_i$ is determined by comparing the orientation of $T_{p_i}(I \times J \times X)$ to that of the sum $d \widehat h(T_{x_i}(I \times J \times K)) \oplus d \widehat h(T_{y_i}(I \times J \times K))$. The \bit{Dax invariant $\Delta(h)$ of the homotopy $h$}, first considered in \cite{dax}, is the sum $$\Delta(h) := \sum_{i=1}^m \epsilon_i g_i \in \mathbb{Z}[\pi_1(X) \setminus 1].$$ Since the Dax invariant of a homotopy is equal to a signed sum of double points, it is clearly additive under the composition of homotopies. For more background, see Gabai \cite[Section 4]{dave:LBL}, as well as Sections $4.1$ and $5.1$ of Kosanovi\'c-Teichner \cite{KT}, including Definition $4.4$. 
\end{definition}

\fig{180}{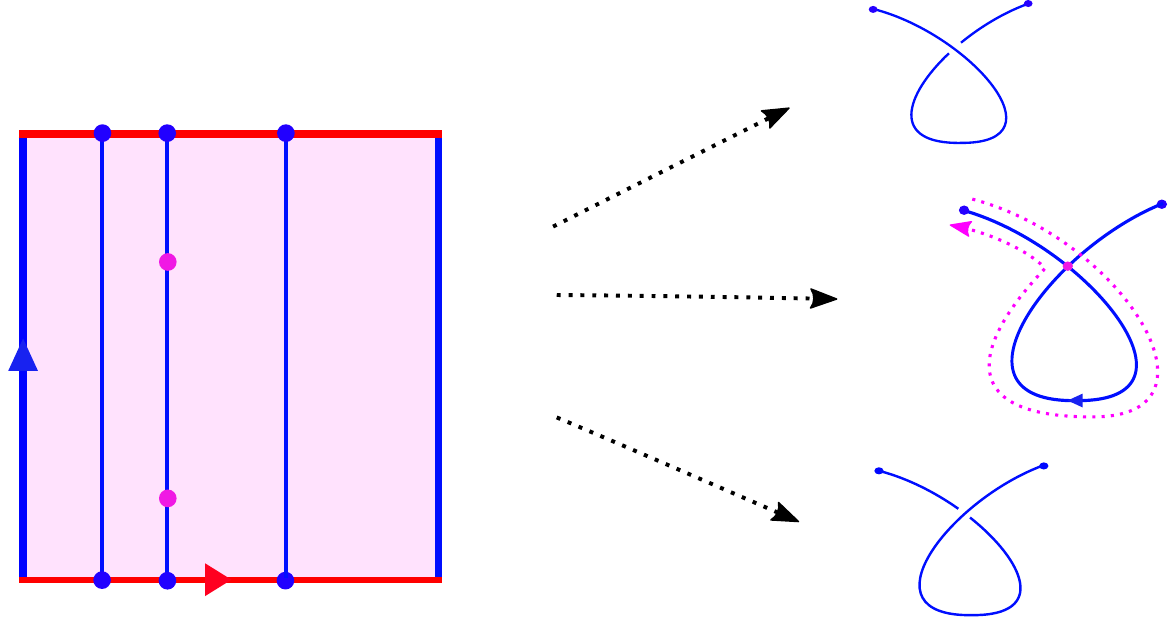}{
\put(-282,-5){\small{$J$}}
\put(-350,72){\small{$K$}}
\put(-285,101){\small{$y_i$}}
\put(-285,32){\small{$x_i$}}
\put(-32,90){\small{$p_i$}}
\put(0,80){\small{$\gamma_i$}}
\put(-143,98){\small{$\widehat h$}}
\caption{For each $t \in I$, the map $\widehat h_t : J \times K \to I \times J \times X$ from Definition \ref{dax} can be thought of as a homotopy in $X$ of the interval $K$, parametrized by $J$. Each double point $p_i$ has a unique associated group element $g_i \in \pi_1(X)$ represented by the based oriented loop $\gamma_i$ pictured on the right of the diagram.}
\label{daxpic}}

Note that having a canonical positive direction along the interval $K$ is critical for the uniqueness of both the sign and group element associated to each double point of the immersion $\widehat h$ above, since this gives an order to the sheets of the $3$-manifold near each intersection point. 

\begin{definition} \label{fqdefn} The \bit{Freedman-Quinn invariant $\fq(h)$ of the homotopy $h$} is the image of $\Delta(h)$ under the quotient map from $\mathbb{Z}[\pi_1(X) \setminus 1]$ to the ring $\mathbb{Z}_2[T]$, where $T \subset \pi_1(X)$ denotes all elements of order $2$. Note that this is not the classical definition from Theorem $10.5$ of \cite{freedman-quinn:4-manifolds}, but rather a new formulation due to Kosanovi\'c-Teichner \cite[Theorem $5.11$]{KT}. For the classical definition of the Freedman-Quinn invariant, refer to the original source \cite[Theorem $10.5A$]{freedman-quinn:4-manifolds}, and for a geometric interpretation, see \cite{me}. This invariant of homotopies, up to deformation, gives rise to the \bit{Freedman-Quinn concordance invariant $\mathrm{fq}(D_0, D_1)$ of a pair of homotopic disks}. In particular, a pair of disks $D_0$ and $D_1$ are related by a homotopy with trivial Freedman-Quinn invariant if and only if the Freedman-Quinn invariant of the pair vanishes.
\end{definition}

As mentioned in the introduction, this concordance invariant was defined (in their context, for pairs of homotopic spheres) by Freedman and Quinn \cite[Theorem $10.5A$]{freedman-quinn:4-manifolds}. Careful expositions can also be found in Section $4$ of Schneiderman-Teichner \cite{st} as well as Section $4$ of Klug-Miller \cite{mandm}. In the presence of an embedded sphere that is \emph{dual} to homotopic spheres $S$ and $T$ (look ahead to Definition \ref{dual}), the Freedman-Quinn invariant $\fq(S,T)=0$ if and only if $S$ and $T$ are smoothly isotopic. This fact is due in combination to Gabai \cite{dave:LBT}, who first showed (albeit without the language of Freedman-Quinn) the ``if" direction, and to Schneiderman-Teichner \cite{st}, who later tied in the Freedman-Quinn invariant to give a complete proof. In contrast, if the common dual is allowed a \emph{non-trivial} normal bundle, it is also required that an additional invariant, the Kervaire-Milnor invariant -- due to Stong \cite{stong} -- vanishes in order to guarantee even a smoothly embedded concordance between the spheres $S$ and $T$. 

\begin{definition} \label{dax2} 
Given homotopic properly embedded disks $D_0$ and $D_1$, the \bit{Dax invariant $\Dax(D_0, D_1)$ of the pair} is an isotopy invariant equal to the image of the Dax invariant of the homotopy $\Delta(h)$ in the quotient $\mathbb{Z}[\pi_1(X) \setminus 1]/ \dax(\pi_3(X))$, where $h$ is any homotopy from $D_0$ to $D_1$, and the map $\dax: \pi_3(X) \to \mathbb{Z}[\pi_1(X) \setminus 1]$ is the \bit{Dax homomorphism} first defined by Dax \cite{dax}. Later, Gabai developed the use of this invariant more explicitly as an isotopy invariant for homotopic disks properly embedded in $4$-manifolds; see Gabai's Definition $3.8$, Theorem $0.3$, and the discussion of the ``Dax kernel"  in \cite{dave:LBL} for instance. The image of $\dax$ consists exactly of the Dax invariant of the self-homotopies of any fixed disk; hence the dependence of the Dax invariant $\Dax(D_0, D_1)$ on the specific homotopy $h$ from $D_0$ to $D_1$ is eliminated in its quotient, as shown in Corollary $0.5$ of Gabai \cite{dave:LBL}. 
\end{definition}

\begin{remark}
The homomorphism $\dax: \pi_3(X) \to \mathbb{Z}[\pi_1(X) \setminus 1]$ has been computed explicitly for few families of $4$-manifolds. Clearly, when $\pi_1(X)=1$ or $\pi_3(X) = 0$, the Dax homomorphism is trivial. On the other hand, the map is known to be non-trivial for some $4$-manifolds with fundamental group $\mathbb{Z}$, in particular $S^1 \times S^3$ and $(S^1 \times S^3) \natural (S^2 \times D^2)$. One should refer to \cite{KT}, \cite{budgab}, \cite{dave:LBL}, and \cite{AS} for greater detail. For related computations of the Dax homomorphism in more general settings, see the work of Kosanovi\'c in \cite{K1} and \cite{K2} for example. 
\end{remark}

\begin{example} \label{Delta} 
Given a properly embedded disk $D \subset X$ with dual $G$ and self-referential disk $D_\omega \subset X$ for some word $\omega \in \mathbb{Z}[\pi_1(X)]$, Theorem $4.9$ of Gabai \cite{dave:LBL} shows that the Dax invariant $\Dax(D, D_\omega)$ is equal to $\mathcal D(\omega)$, where $\mathcal D: \mathbb{Z}[\pi_1(X)] \to \mathbb{Z}[\pi_1(X) \setminus 1]/\dax(\pi_3(X))$ is the group homomorphism induced by sending $1 \mapsto 0$ and $g \mapsto g + \bar g$  in the quotient\footnote{We use the notation $\bar g := g^{-1}$ (without a sign change).} for each non-trivial $g \in \pi_1(X)$. This can also be seen from the group action given in Kosanovi\'c and Teichner \cite[Theorem A]{KT}. Although again it is relevant to note that Gabai is working with a dual sphere in the boundary of the $4$-manifold, his proof of Theorem $4.9$ can be carried out in our present setting where the dual lies in the interior of $X$. In fact, the homotopy depicted in Figures $11-15$ of Gabai in \cite{dave:LBL} is exactly the homotopy illustrated in our Figure \ref{Dg} whose Dax invariant we compute in Example \ref{exampledax}.
\end{example}

\section{Modifying regular homotopies}

We apply techniques due to Freedman and Quinn in \cite{freedman-quinn:4-manifolds}, and Quinn in \cite{quinn}, that can be used in the presence of a common dual to ``improve" a regular homotopy $h: I \times D^2 \to X$ by deforming it via a homotopy rel boundary. We begin by proving the following analog, for properly embedded disks, of Lemma $6.1$ of Gabai \cite{dave:LBT} and Lemma $2.1$ of Schneiderman-Teichner \cite{st}. 

\fig{130}{fqmove3.png}{
\put(-324,115){\small{$h(\gamma)$}}
\put(-189,72){\small{$h(\alpha)$}}
\put(-182,120){$c$}
\put(-260,15){$A$}
\put(-285,31){$\alpha$}
\put(-128,13){$U$}
\put(-128,37){$K$}
\put(-225,32){w.m.}
\put(-412,72){$G$}
\put(-216,111){$\Delta$}
\put(-400,96){\small{Im$(h)$}}
\caption{The (immersed) Whitney move used in our proof of Lemma \ref{homotopyrelG} to modify the intersection between the image of $h$ and the dual sphere $G$, shown from the perspective of $X$ (top diagram) as well as from the pre-image $I \times D^2$ of the homotopy (bottom diagram). Performing a Whitney homotopy of $\mathrm{Im}(h)$ along the immersed disk $\Delta \subset X$ deforms the homotopy $h$ so that the arc $A$ in the pre-image of the double curves in $I \times D^2$ is replaced by an unknotted arc $U$ and a knot $K$, the result of surgering along the arc $\alpha$. }
\label{fq2pic}}

\begin{lemma} \label{homotopyrelG} 
Let $h$ be a homotopy from one properly embedded disk $D_0$ to another $D_1$ properly embedded in $X$. If $D_0$ and $D_1$ share a common $\pi_1$-negligible dual $G$ in the interior of $X$, then there is a homotopy $h'$ supported away from $G$ and $\partial X$ with $\mathrm{fq}(h)=\mathrm{fq}(h')$. 
\end{lemma}

\proof 
After a perturbation of the homotopy $h: I \times D^2 \to X$, its image in $X$ intersects the dual $G$ in a collection of embedded curves and a single arc. The pre-image of these intersections is the union of a link $L$ and a properly embedded arc $A$ in $I \times D^2$. We shall abuse notation below by continuing to refer to the link as $L$ and the arc as $A$ even after several modifications to the original homotopy.

First, we deform the homotopy $h$ in order to unknot the arc $A$. Start by choosing an arc $\alpha$ embedded in $I \times D^2$ with boundary points $x, y$ on $A$ but interior disjoint from $A \cup L$. Choose $\alpha$ so that surgering $A$ along $\alpha$ results in both an unknotted arc $U$ and a knot $K \subset I \times D^2$, as illustrated in the bottom diagram of Figure \ref{fq2pic}. Let $c \subset G$ be an embedded arc running parallel to the arc $h(A)$ and connecting the pair of points $h(x)$ and $h(y)$. 

The union of the arcs $c$ and $h(\alpha)$ give a closed, embedded loop in $X$, shown in the top diagram of Figure \ref{fq2pic}. This loop is null-homotopic in $X$: this can be shown using the facts that the points $x$ and $y$ lie on the pre-image of the same double curve of $h$, and both $G$ and $I \times D^2$ are simply-connected. Therefore, the loop $c \cup h(\alpha)$ bounds an immersed disk $\Delta \subset X$. Performing ``boundary twists" of the disk $\Delta$ around $G$, as in \cite[\S1.3--1.4]{freedman-quinn:4-manifolds}, arranges for a neighborhood of $\Delta$ to have the local model shown in Figure \ref{fq2pic} but with the interior of $\Delta$ immersed. Each boundary twist adds an extra point of intersection to $\mathrm{int}(\Delta) \pitchfork G$, but this can be removed by tubing $\Delta$ over an immersed dual $G'$ to the sphere $G$. Recall, such a dual exists since $G$ is $\pi_1$- negligible. In fact, the disk $\Delta$ may be tubed repeatedly over $G'$ until it is disjoint from $G$. 

\fig{160}{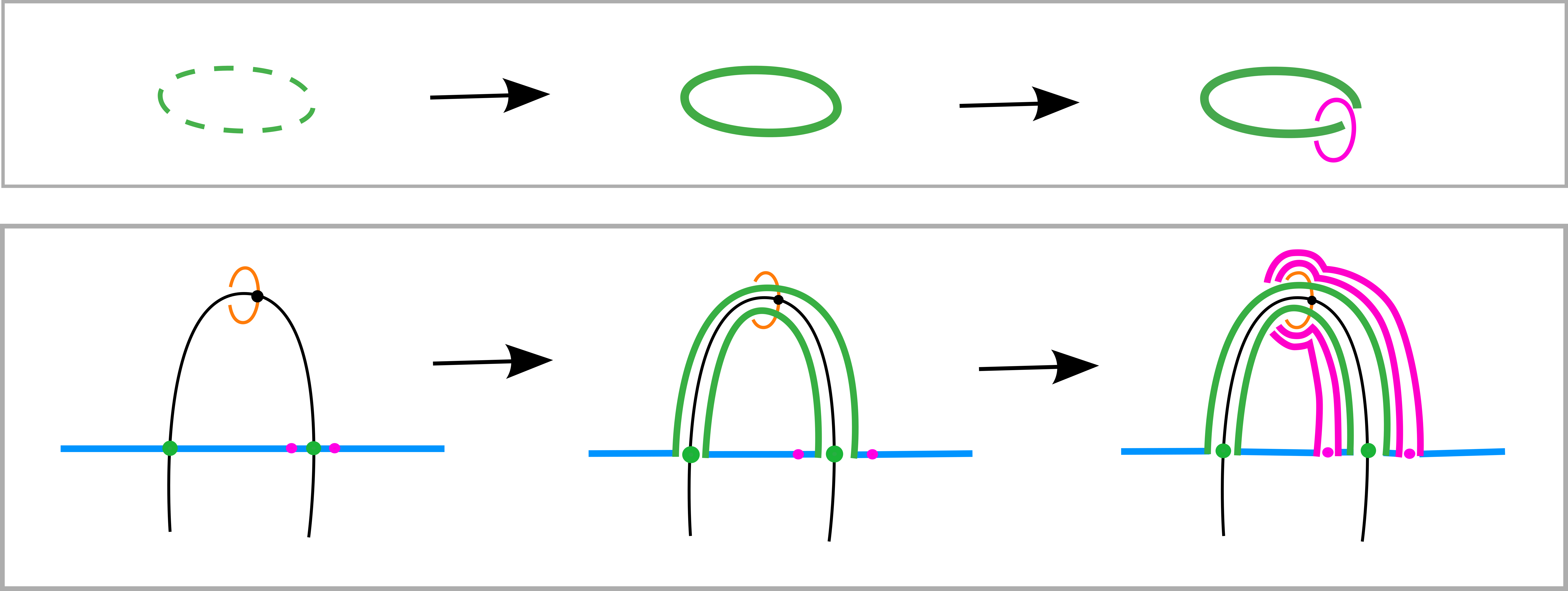}{
\put(-380,119){$\gamma$}
\put(-243,118){$0$}
\put(-100,117){$0$}
\put(-55,116){$0$}
\put(-50,10){$G$}
\put(-195,11){$G$}
\put(-336,12){$G$}
\put(-410,26){\small{Im$(h)$}}
\put(-125,26){\small{Im$(h')$}}
\put(-352,82){\small{$G'$}}
\caption{The \emph{surgery operation} used in Lemma \ref{homotopyrelG} to replace a homotopy $h$ with one $h'$ whose image has one less closed curve of intersection with the dual sphere $G$. For each component $\gamma$ of $L$, first perform ambient $0$-framed surgery on the immersed $3$-manifold $\mathrm{Im}(h)$ as in the bottom diagram, replacing a neighborhood of $h(\gamma)$ with a solid torus (green) whose meridianal disk is parallel to one in $G$ bounded by $h(\gamma)$. The top diagram illustrates the effect of this surgery on the (possibly non-trivial) curve $\gamma$ in the pre-image. To ``undo" this change to the domain, ambiently surger the immersion in $X$ along the image of the meridian to $\gamma$, replacing its tubular neighborhood by a solid torus (pink) whose meridianal disk passes over the immersed sphere $G'$ as depicted in the bottom diagram, to avoid hitting the dual $G$. This ambient surgery gives an immersion of $I \times D^2$, and so this operation produces a second homotopy $h'$.}
\label{fq1pic}}

An immersed Whitney homotopy pushing the image of $h$ across the immersed disk $\Delta$ modifies the curves in $\mathrm{Im}(h) \pitchfork G$ so that the arc $A$ in the domain $I \times D^2$ is surgered along $\alpha$. A similar construction, in a slightly different context, was used by Freedman and Quinn \cite[pg. $189$]{freedman-quinn:4-manifolds}. This deformation of $h$ is essentially a generalization of the procedure of ``summing" that we use later, shown in Figure \ref{steps}.  By construction, performing this immersed Whitney move along $\Delta$ unknots the arc $A$, and also adds one knotted component to the link $L \subset I \times D^2$. Indeed, after a reparametrization of the homotopy $h$, we may assume that $\alpha = I \times \{\pt\}$. Note that the homotopy at this stage is simply a deformation of the original one, and so its Freedman-Quinn invariant is unchanged. 

To remove all closed curves of intersection between the dual $G$ and the image of $h$, we perform {surgery replacements} to $h$ for each component $\gamma \subset L$, as shown in Figure \ref{fq1pic}. When no such curves remain, this gives a homotopy $h'$ supported away from $G$, as desired. We note again that a similar construction was used in a slightly different context by Freedman and Quinn \cite{freedman-quinn:4-manifolds}; compare their illustration on page $191$ to ours. 

The double curves of this new homotopy $h'$ include those of $h$, as well as two new sets of double curves. First, near each self-intersection of $G'$, the immersed $2$-handles of the image of $h'$ from Figure \ref{fq1pic} intersect both themselves and each other, since their cores run parallel to $G'$. In addition, the image of $h'$ intersects itself in curves parallel to those where the sphere $G'$ intersects the original homotopy $h$; there is one pair of such curves for each $2$-handle of $\mathrm{Im}(h')$ running along $G'$.  

In either of these two cases, the pre-image of each double curve of $h'$ that is not already a double curve of $h$ has \emph{two} connected components in $I \times D^2$. Therefore, the lift $\widehat {h'}$ from Definition \ref{dax} must have an \emph{even} number of double points corresponding to each non-trivial element of $\pi_1(X)$, besides those it inherits from $\widehat{h}$. This implies that the Freedman-Quinn invariant of the homotopy $h'$ is equal to that of $h$, by Definition \ref{fqdefn}. 
 \qed 
\smallskip

Mainly for our use later in our proof of Theorem \ref{one}, we are interested in improving a regular homotopy so that its double points can be organized in the following sense. 

\begin{definition}\label{sheets}
Let $\Sigma \subset X$ be a smoothly immersed disk, whose boundary is embedded in the boundary of $X$. The pre-images $x_i, x_i' \in D^2$ of  each of the double points $p_1, \dots, p_n$ of $\Sigma$ can be partitioned arbitrarily into two sets $\{x_1,  \dots, x_n\}$ and $\{x_1', \dots, x_n'\}$ each of which is called a \bit{sheet} of the immersion. See Figure \ref{sheetspic}, and compare the related discussion of ``sheets" in \cite{st}.  If $\Sigma \subset X$ is the middle level of a homotopy $h$ consisting of $n$ finger and $n$ Whitney moves, then it has positive double points $p_1, \dots, p_n$ and negative double points $q_1, \dots, q_n$ with pre-images  $x_i, x_i' \in D^2$ and $y_i, y_i' \in D^2$, respectively. A choice of sheets $\{x_1, y_1, \dots, x_n, y_n\}$ and $\{x_1', y_1', \dots, x_n', y_n'\}$ such that each of the pairs of points $\{x_i, y_i\}$ and $\{x_i', y_i'\}$ are connected by the pre-image of a descending (resp. ascending) Whitney arc are called \bit{descending (resp. ascending) sheets}. 
\end{definition}

\fig{200}{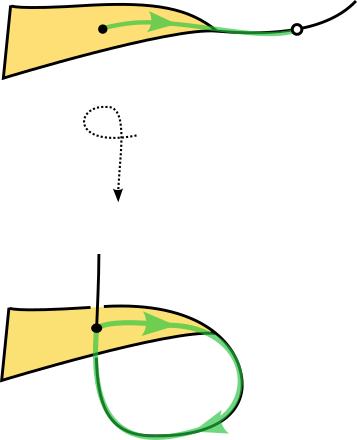}{
\put(-178,180){$D^2$}
\put(-174,40){$\Sigma$}
\put(-130,47){$p$}
\put(-128,183){$x$}
\put(-30,193){$x'$}
\put(-44,25){$\widehat \gamma$}
\caption{A double point $p$ of an immersed disk $\Sigma$ with pre-images $x, x' \in D^2$ corresponds to an element $g \in \pi_1(X)$ represented by the based loop given by dragging any point on the oriented curve $\widehat \gamma$ to the basepoint of $\pi_1(X)$. In general, since there is no canonical way to label the pre-images $x$ and $x'$, the double point $p$ corresponds to the element $\bar g \in \pi_1(X)$ as well. However, a choice and ordering of a pair of sheets for the pre-images of the double points of $\Sigma$ determines a \emph{unique} group element corresponding to each double point $p$, since in this case $\widehat \gamma$ has a well-defined orientation as the image of an oriented path in $D^2$ from the pre-image $x$ in the first sheet to the pre-image $x'$ in the second sheet.}
\label{sheetspic}}

Recall from Section \ref{daxsection} that the basepoint $b$ of $\pi_1(X)$ is assumed to lie in $\partial X$ on the boundary circle of the properly immersed disk. In this case, each double point $p$ of the middle level $\Sigma$, with pre-images $x$ and $x'$, \emph{corresponds} to a group element which we denote $\gamma(p) \in \pi_1(X)$, represented by the closed oriented loop in $X$ formed from the image of two paths in the disk $D^2$: one running from the pre-image of the basepoint $b$ to $x$, and the other from $x'$ to the pre-image of $b$. Since the disk $D^2$ is $1$-connected, all such loops are homotopic, up to reversal of orientation. In general there is no canonical ordering of the pre-images $x$ and $x'$, and so the group element $\gamma(p)$ is well-defined only up to its inverse. However, by specifying a choice of sheets as in Definition \ref{sheets}, all double points can be assigned a group element uniquely by fixing an ordering of the sheets. This is more carefully explained in the caption below Figure \ref{sheetspic}.

\begin{definition}
Consider sheets $\{x_1, y_1, \dots, x_n, y_n\}$ and $\{x_1', y_1', \dots, x_n', y_n'\}$ of an immersed disk $\Sigma \subset X$ as in Definition \ref{sheets}. A pair of sheets $\{x_1, y_1, \dots, \ast_i', \dots, x_n, y_n\}$ and $\{x_1', y_1', \dots, \ast_i, \dots, x_n', y_n'\}$ where $\ast_i \in \{x_i, y_i\}$  is said to be obtained from the first pair of sheets by a \bit{swap}. Each swap is assigned both the sign and corresponding group element of the double point $h(\ast_i)=h(\ast_i')$. From now on, we assume implicitly that the both pairs of sheets are ordered so that this group element is well-defined.
\end{definition} 

The following is a direct consequence of Schneiderman-Teichner \cite[Lemma $4.6$]{st} (stated in the source for spheres, but their proof extends to homotopies of disks). We state it here as a lemma so that it can be more easily referred to later, using the terminology from this section. 

\begin{lemma}[\cite{st}]\label{stlemma}
If $\fq(h)=0$, then any ascending and descending sheets of $h$ are related by an equal number of positive and negative swaps corresponding to each element of $\pi_1(X)$.
\end{lemma} 

\fig{180}{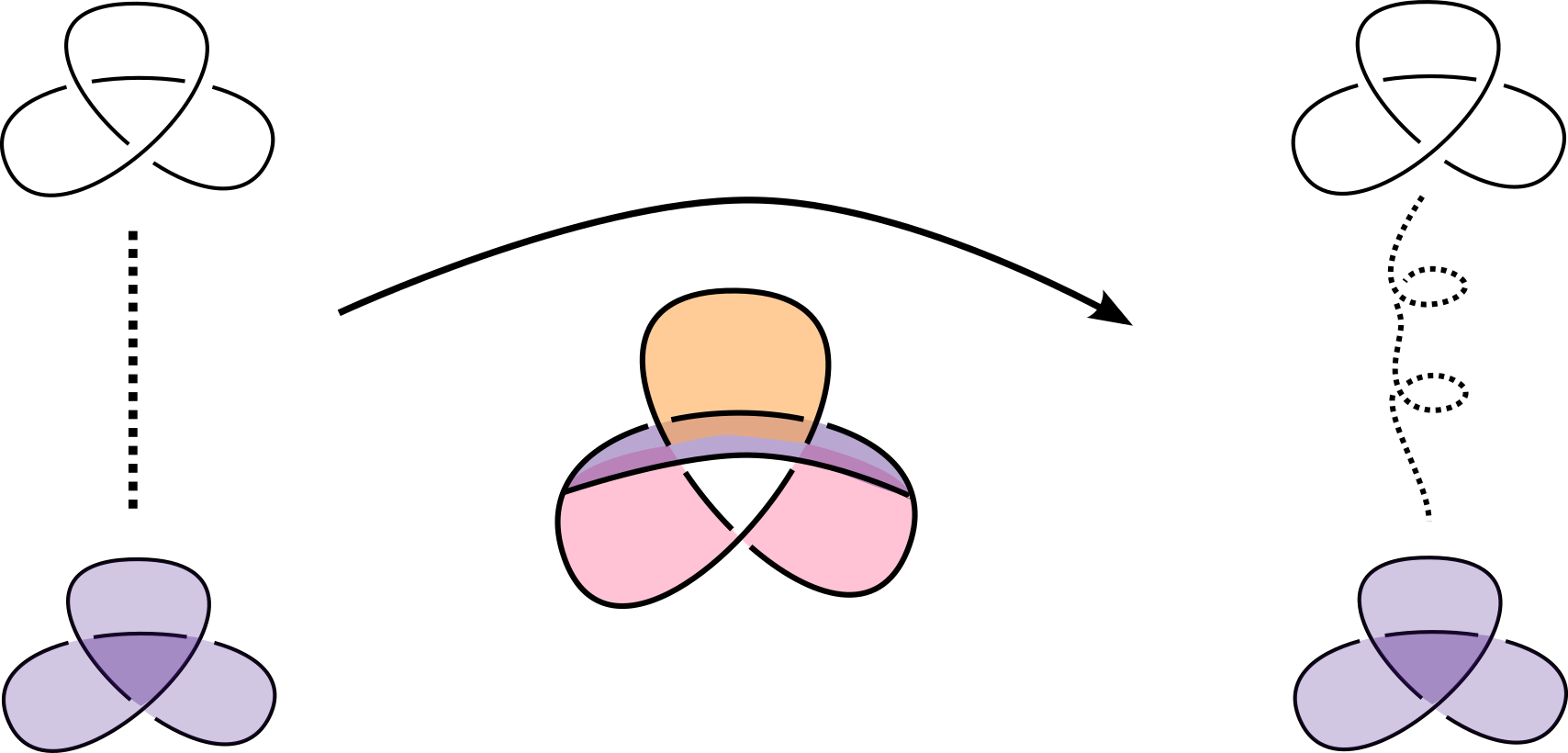}{
\put(-227,147){$1$ f.m. + $1$ w.m.}
\put(-335,95){$\times ~I$}
\put(-20,95){$2$ ``twists"}
\caption{A regular homotopy between trivials disks properly embedded in the $4$-ball, consisting of one finger move (f.m.) and one Whitney move (w.m.). One can check that the ascending and descending Whitney disks shaded in the picture of the middle level $\Sigma$ give ascending and descending sheets of double points related by a single swap corresponding to the trivial element. By \cite{me}, this homotopy takes the trivial disk on the left to the trivial disk on the right that differs only by inserting two full rotations of the unknot in each level set of the $I$ factor. However, this does not effect the isotopy class of the resulting disk rel boundary in $B^4$. For more detailed illustrations of such a homotopy, see Figures $14$-$16$ of \cite{me}.}
\label{hompic}}

Although Schneiderman-Teichner \cite{st} state this fact explicitly only for \emph{non-trivial} elements of $\pi_1(X)$, any homotopy can be deformed within a $4$-ball as in Figure \ref{hompic} so that there are an equal number of positive and negative swaps corresponding to the identity as well.  We are interested in proving the related lemma:

\begin{lemma}\label{sheetlemma}
Suppose a homotopy $h: I \times D^2 \to X$ with $\fq(h)=0$ is supported away from some $\pi_1$-negligible sphere $G \subset X$ dual to the disks $h_t(D)$ for each $t \in I$.  Then, $h$ can be deformed through homotopies (rel boundary and away from $G$) to a homotopy $h'$ whose middle level admits a pair of sheets that is both \emph{ascending} and \emph{descending}, as in Definition \ref{sheets}.  
\end{lemma} 

In order to obtain this result, we will need an operation on Whitney disks defined by Quinn in \cite[pg. 355]{quinn} that can be used to physically realize pairs of sheet swaps corresponding to the same group element, but of opposite sign. Note that it is crucial in this construction that the homotopy $h$ is supported away from the dual sphere $G$. 

\begin{definition}\label{defsum}
Let $W_1$ and $W_2$ be Whitney disks of the homotopy $h$, with Whitney arcs $\omega_1, \omega_1^*$ and $\omega_2, \omega_2^*$. Let $\sigma$ be a loop in $X$ given by the concatenation of embedded paths: {\bf(i)} along $\Sigma$ from the interior of $\omega_1$ to the interior of $\omega_2$, {\bf (ii)} along $W_2$ from the interior of $\omega_2$ to the interior of $\omega_2^*$, {\bf (iii)} along $\Sigma$ from the interior of $\omega_2^*$ to the interior of $\omega_1^*$, and finally {\bf(iv)} along $W_1$ from the interior of $\omega_1^*$ to the point in the interior of $\omega_1$ where path (i) began. If in addition the interior of each path is disjoint from all Whitney arcs of $h$, then we call the loop $\sigma$ a \bit{square} for these Whitney disks. Suppose such a loop bounds a disk $S$ with local model as shown in Figure \ref{summing}, embedded away from the dual $G$, the middle level $\Sigma$, as well as the Whitney disks of $h$. The homotopy $h'$ obtained by \bit{summing} the Whitney disks $W_1$ and $W_2$ is the homotopy consisting of exactly the same finger moves as $h$, and Whitney disks equal to those of $h$ except with the disks $W_1'$ and $W_2'$ from Figure \ref{summing} replacing $W_1$ and $W_2$.
\end{definition}

This summing construction is critical to the proof of Lemma \ref{sheetlemma} to ensure that changes of doublepoint pairings can be realized by a deformation through homotopies. Quinn \cite[pg. 355]{quinn} explicitly deforms the homotopy $h$ to the homotopy $h'$ rel boundary\footnote{Although Quinn works in a simply-connected setting to find the disk $S$, the details of his deformation from $h$ to $h'$ do not rely on the fundamental group of the ambient manifold.}. The same deformation can also be thought of as a Whitney homotopy in 5-dimensions between the \emph{traces} of the homotopies in $I \times X$ as in Freedman and Quinn \cite[Section $10.9$, pp. 189--190]{freedman-quinn:4-manifolds}, or from the perspective of the figure at the bottom of \cite[pp. 182]{freedman-quinn:4-manifolds} in which the existence of the disk $C_1$ is essentially the same as the existence of our disk $S$. With this operation in mind, we proceed with the proof of the lemma. 

\smallskip
\noindent
\emph{Proof of Lemma \ref{sheetlemma}.} Let  $\{x_1, y_1, \dots, x_n, y_n\}$ and $\{x_1', y_1', \dots, x_n', y_n'\}$ be ascending sheets of the homotopy $h$. Since $\fq(h)=0$, by Schneiderman and Teichner's Lemma \ref{stlemma} above, these sheets are related to any choice of descending sheets for $h$ by an equal number of positive and negative swaps corresponding to each element of $\pi_1(X)$. We proceed by induction on the number $n$ of cancelling pairs of swaps needed. 

If no swaps are required, then we are done. So, assume that there are $n>0$ cancelling pairs of swaps needed to convert the given ascending sheets to the descending sheets. First, suppose one such pair of swaps exchanges the sheets of the pre-images $x_i, x_i'$  and $y_i, y_i'$. Since we started with a set of ascending sheets for $h$, both of the pairs of points $\{x_i, y_i\}$ and $\{x_i', y_i'\}$ are connected by the pre-image of an ascending Whitney arc. So, a new set of ascending sheets for $h$ is given by $\{x_1, y_1, \dots, x_i', y_i', \dots, x_n, y_n\}$ and $\{x_1', y_1', \dots, x_i, y_i, \dots, x_n', y_n'\}$. These ascending sheets are related to the fixed set of descending sheets by $n-1$ cancelling pairs of swaps. 

Now, suppose that \emph{all} cancelling pairs of swaps exchange the sheets of the pre-images $x_i, x_i'$ and $y_j, y_j'$ for $i \not = j$. Let $\sigma$ be a square for the Whitney disks $W_i$ and $W_j$, as in Definition \ref{defsum}, one of whose sides connects the Whitney arcs with boundaries $\{x_i, y_i\}$ and $\{x_j', y_j'\}$, and another that connects the Whitney arcs with boundaries $\{x_i', y_i'\}$ and $\{x_j, y_j\}$. The loop $\sigma$ is null-homotopic in $X$ since the cancelling pair of swaps corresponds to the same element of $\pi_1(X)$, and in fact, bounds an immersed disk in $X-G$ since $G$ is $\pi_1$-negligible. Such a disk can be made to have the framing needed to admit the local model in Figure \ref{summing} by ``boundary twisting" as in \cite[\S1.3--1.4]{freedman-quinn:4-manifolds} at the cost of creating more intersection points with the middle level $\Sigma$ of the homotopy. 

\fig{160}{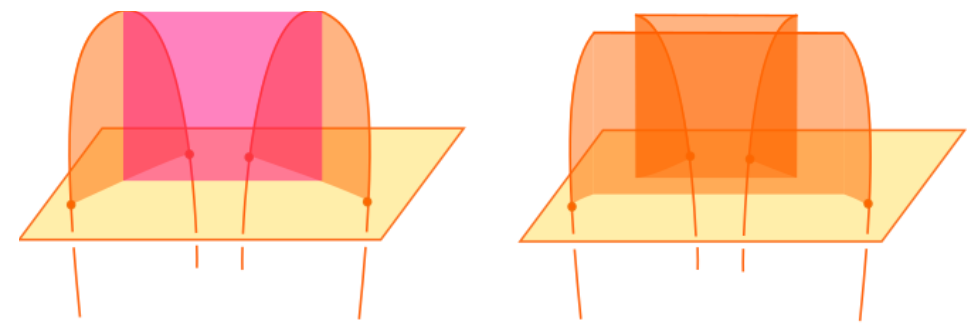}{
\put(-80,110){$W_1'$}
\put(-130,125){$W_2'$}
\put(-365,130){$S$}
\put(-427,113){$W_1$}
\put(-310,113){$W_2$}
\put(-284,32){$\Sigma$}
\put(-46,30){$\Sigma$}
\caption{The local model, originally given by Quinn \cite[pg. 355]{quinn}, needed to sum two Whitney disks $W_1$ and $W_2$ to produce new Whitney disks $W_1'$ and $W_2'$. The embedded disk $S$ along which the sum is performed has boundary consisting of two boundary arcs on $\Sigma$, one on $W_1$, and one on $W_2$. The disk $S$ is also embedded away from the dual $G$, the middle level $\Sigma$, and the Whitney disks of $h$. Its normal disk bundle restricts to the framing on the boundary indicated by the model, so that $W_1'$ and $W_2'$ can be formed from the union of halves of $W_1$ and $W_2$ with disjoint parallel push-offs of $S$.}
\label{summing}}

Now, we homotop the resulting immersed disk $D$ to an embedded disk disjoint from both the Whitney disks of the homotopy $h$ and its middle level $\Sigma$, by performing the sequence of regular homotopies shown in Figure \ref{seqofhom}. Summing the Whitney disks $W_i$ and $W_j$ using the resulting disk we call $S$, with $\partial S= \sigma$, gives new Whitney disks $W_i'$ and $W_j'$ whose Whitney arcs have boundaries $\{x_i, y_j'\}$, $\{x_i', y_j\}$, $\{x_j, y_i'\}$, and $\{x_j', y_i\}$ in the pre-image. Therefore, the homotopy $h'$ obtained by summing the Whitney disks $W_i$ and $W_j$ has ascending sheets given by  $\{x_1, y_1, \dots, x_i', \dots, y_j', \dots, x_n, y_n\}$ and $\{x_1', y_1', \dots, x_i, \dots, y_j, \dots, x_n', y_n'\}$. As in the first case, this set of ascending sheets is related to the fixed set of descending sheets by $n-1$ cancelling pairs of swaps. Hence, our inductive argument is complete. \qed 

\fig{150}{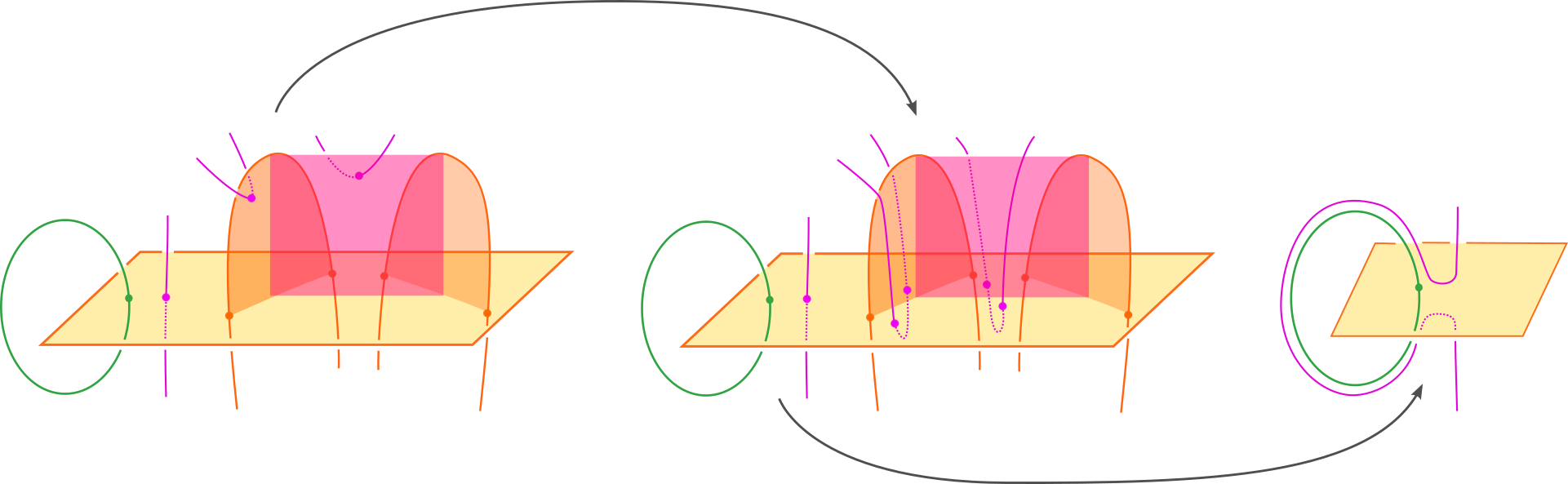}{
\put(-380,80){$D$}
\put(-470,85){$G$}
\put(-328,156){Steps $1$ and $2$}
\put(-320,50){$\Sigma$}
\put(-155,-11){Step $3$}
\caption{Starting with immersed disk $D$, the following steps create the embedded disk $S$ needed for the local model in Figure \ref{summing}. {\bf Step $1$:} Homotop the disk $D$ to an embedding, by pushing one sheet of each of its double points across its boundary and through the surface $\Sigma$. {\bf Step $2$:} Convert all intersections between $D$ and the Whitney disks of $h$ into intersections between $D$ and $\Sigma$, by again pushing $D$ across the boundary of the Whitney disks and through $\Sigma$. {\bf Step $3$:} Eliminate all intersection points in $D \pitchfork \Sigma$ by tubing $D$ over the embedded dual $G$. Since both $D$ and the Whitney disks of the homotopy are disjoint from $G$, this does not create any new self-intersections or double points.}
\label{seqofhom}}

\section{Self-referential disks with a common dual}  \label{2}

The focus of this section is the proof of the lemma below, which not only answers Question $5.7$ posed by Gabai in \cite{dave:LBL}, but is also critical to the proof of Theorem \ref{one}, our 4D light bulb theorem for disks.  Recall the definition of a self-referential disk $D_\omega$ from Definition \ref{selfref} and Figure \ref{disks}. 

\begin{lemma}\label{Dgiso}
Let $D$ be a properly embedded disk in $X$. If $D$ has a dual $G$, then for any $g \in \pi_1(X)$ there is a smooth isotopy rel boundary between the self-referential disks $D_g$ and $D_{{\bar{g}}}$ supported away from their common dual $G$. 
\end{lemma}

\proof The finger moves of the regular homotopies $h_g$ and $h_{\bar g}$ given in Figure \ref{Dg} take $D_g$ and $D_{\bar g}$ to the same immersed disk $\Sigma$. Therefore, the descending shadows (as in Definition \ref{shadow}) of both homotopies are tubed disks $(\Sigma, \alpha_g)$ and $(\Sigma, \alpha_{\bar g})$ whose tubing arcs $\alpha_g$ and $\alpha_{\bar g}$ are pictured in the bottom leftmost and rightmost diagrams of Figure \ref{steps}. As noted below Definition \ref{shadow}, it is due to Gabai \cite{dave:LBT} that the descending shadows are smoothly isotopic to the self-referential disks $D_g$ and $D_ {\bar g}$, respectively. Therefore, Figure \ref{steps} gives an isotopy between $D_g$ and $D_ {\bar g}$ viewed as tubed disks, using the tube isotopies A and B from Definition \ref{AB}. Note that each tube isotopy is, indeed, supported away from the dual $G$.
\qed

\fig{165}{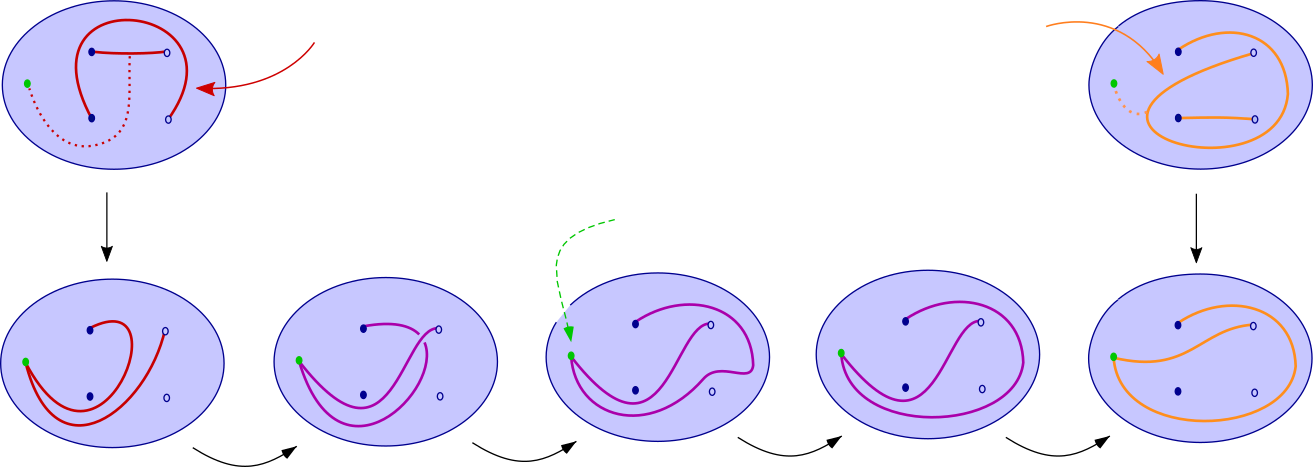}{
\put(-247,94){\small{Pre-image of}}
\put(-243,84){\small{$z=D \pitchfork G$}}
\put(-35,-7){\small{Descending}}
\put(-37,-17){\small{shadow of $h_{\bar g}$}}
\put(-470,-10){\small{Descending}}
\put(-472,-20){\small{shadow of $h_g$}}
\put(-422,83){\small{Tube over $z$}}
\put(-100,83){\small{Tube over $z$}}
\put(-355,157){\small{Pre-image of descending}}
\put(-350,147){\small{Whitney arcs for $h_g$}}
\put(-190,162){\small{Pre-image of descending}}
\put(-185,152){\small{Whitney arcs for $h_{\bar g}$}}
\put(-398,-10){\small{Re-order}}
\put(-300,-8){\small{Isotopy B}}
\put(-205,-7){\small{Isotopy A}}
\put(-110,-7){\small{Isotopy A}}
\caption{A tube isotopy between the descending shadows of $h_g$ and $h_{\bar g}$, defined through a sequence of tube diagrams as in Figure \ref{moves2}. The regular homotopies  $h_g$ and $h_{\bar g}$, as well as the descending Whitney arcs used in the construction of the tubing arcs for their descending shadows, are depicted in Figure \ref{Dg}.}
\label{steps}}

\begin{corollary} \label{middlecor}
Let $D$ be a properly embedded disk in $X$ with dual $G$. Then, for any $\omega \in \mathbb{Z}[\pi_1(X)]$ such that $\mathcal{D} (\omega)= 0$\footnote{The map $\mathcal{D}$ is defined in Example \ref{Delta}.}, there is a smooth isotopy between $D$ and $D_\omega$ supported away from $G$. 
\end{corollary}

\proof
By definition, the condition $\mathcal D(\omega)= 0$ implies that $\omega = \sum \epsilon_i g_i$ such that $\sum \epsilon_i (g_i + \bar g_i)=0$ for $\epsilon_i = \pm 1$ and $g_i \in \pi_1(X)$. Therefore, the terms in the sum $\omega$ can be partitioned into oppositely signed pairs $\{g_i, g_j\}$ with either $g_i= g_j$ or $g_i = \bar g_j$. Lemma \ref{Dgiso} above applied ``locally" gives that the self-referential disk $D_{g - \bar g}$ is isotopic to the self-referential disk $D_{g-g}$ for each $g \in \mathbb{Z}[\pi_1(X)]$. However, by Gabai \cite[Lemma $2.15$]{dave:LBL} the isotopy class of a self-referential disk $D_\omega$ is well-defined for any $\omega \in \mathbb{Z}[\pi_1(X)]$. So, each $D_{g-g}$ is isotopic to $D$, from which it follows that the disk $D_\omega$ is isotopic to $D$ as well. \qed

\section{Converting homotopies to isotopies} \label{3}

As in the previous sections, let $D_0$ and $D_1$ denote disks with a $\pi_1$-negligible common dual sphere $G$, properly embedded in a compact $4$-manifold $X$, that are homotopic rel boundary. We are now ready to prove our main result. 

\begin{theorem} \label{one}
If $\mathrm{fq}(D_0, D_1)=0$, then the disk $D_1$ is smoothly isotopic, rel boundary and away from the dual $G$, to a self-referential disk $(D_0)_\omega$ for some word $\omega \in \mathbb{Z}[\pi_1(X)]$ such that $\mathcal D(\omega) = \rm{Dax}(D_0, D_1)$ in the quotient $\mathbb{Z}[\pi_1(X)]/ \rm{dax}(\pi_3(X))$.
\end{theorem}

\pf 
First note that there exists a homotopy $h$ from $D_0$ to $D_1$ with $\fq(h)=0$ and support away from the dual $G$, by Lemma \ref{homotopyrelG}. After a deformation rel boundary and away from $G$, we may also assume, by Lemma \ref{sheetlemma}, that the pre-images $x_1, x_1', \dots, x_n, x_n'$ and $y_1, y_1' \dots, y_n, y_n'$ of the $n$ positive and $n$ negative double points of the middle level $\Sigma$ of $h$ can be partitioned into two sheets $\{x_1, y_1, \dots, x_n, y_n\}$ and $\{x_1', y_1', \dots, x_n', y_n'\}$ that are both ascending and descending sheets, as in Definition \ref{sheets}. 

\begin{figure}[ht]
    \begin{subfigure}{0.475\textwidth}
      \centering
      \begin{overpic}[width=0.95\textwidth]{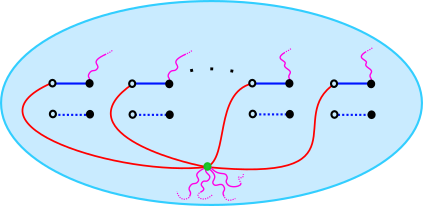}
      \end{overpic}
      
      \caption{Tubing arcs $\alpha_1, \dots, \alpha_n$ (red) drawn in ``standard position" with respect to the descending Whitney arcs $\omega_1, \dots, \omega_n$ (blue solid) and $\omega_1', \dots, \omega_n'$ (blue dashed). The other half of the tubing arcs  $\beta_1, \dots, \beta_n$ (purple) are disjointly embedded from the $\alpha_i$ but not necessarily the Whitney arcs. 
      \label{fig1}}
    \end{subfigure}
    \hfill
    \begin{subfigure}{0.475\textwidth}
      \centering
      \begin{overpic}[width=0.95\textwidth]{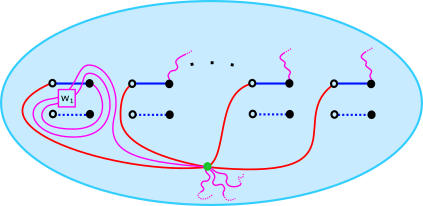}
      \end{overpic}
      
      \caption{A sequence of tube isotopies moving only the arc $\beta_1$ arrange it into the ``standard position" illustrated above, disjoint from all Whitney arcs $\omega_i, \omega_i' \not = \omega_1$ and crossing underneath all arcs $\beta_i$ with $i>1$. The index $w_1 \in \mathbb{Z}$ denotes the winding number of the (oriented) loop $\alpha_1  \omega_1 \beta_1$ based at the pre-image of $\Sigma \pitchfork G$ around the arc $\omega_1'$. \label{fig2}}
    \end{subfigure}
    \vskip\baselineskip
    \begin{subfigure}{0.475\textwidth}
      \centering
      \begin{overpic}[width=0.95\textwidth]{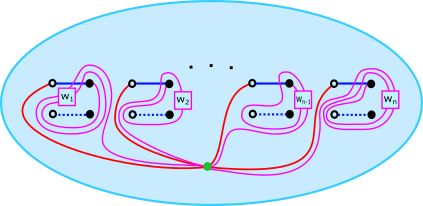}
      \end{overpic}
      
      \caption{A sequence of tube isotopies analogous to the one from Figure \ref{fig2} moves the remaining arcs $\beta_2, \dots, \beta_n$ into ``standard position" as illustrated above. As with $i=1$, the index $w_i \in \mathbb{Z}$ denotes the winding number of the oriented loop $\alpha_i \omega_i \beta_i$ based at the pre-image of $\Sigma \pitchfork G$ around the arc $\omega_i'$.
      \label{fig3}}
    \end{subfigure}
    \hfill
    \begin{subfigure}{0.475\textwidth}
      \centering
      \begin{overpic}[width=0.95\textwidth]{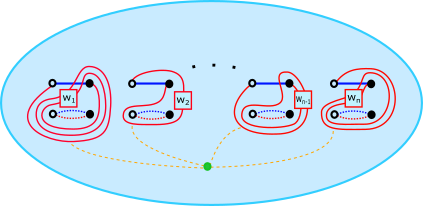}
      \end{overpic}
      
      \caption{Pre-images of the ascending and descending Whitney arcs for the homotopy shown in Figure \ref{twistdisks}. The tube diagram from Figure \ref{fig3} is one for the descending shadow of the homotopy -- this can be seen by following the construction from Definition \ref{shadow}, using the dashed orange arcs as the $\nu_i$.
      \label{fig4}}
    \end{subfigure}
    
\caption{
    Illustrating the proof of Theorem \ref{one}.
    }
    \label{4x4fig}
\end{figure}

Let $\mathcal A$ be an ascending shadow, as in Definition \ref{shadow}, of the homotopy $h$ after this deformation. As noted below Definition \ref{shadow},  the disk $\mathcal A$ is isotopic to $D_1$ by Gabai \cite[Lemma 5.1]{dave:LBT}. In particular, because the sheets above are both ascending and descending, the disk $\mathcal A$ is equal to an embedded tubed disk $(\Sigma, \alpha)$ with disjointly embedded tubing arcs $\alpha=\{\alpha_1, \beta_1, \dots, \alpha_n, \beta_n\}$ that can be chosen so that their pre-images are in the ``standard position" shown in Figure \ref{fig1}, with each $\alpha_i$ connecting $x_i$ to the pre-image of $\Sigma \pitchfork G$, and each $\beta_i$ connecting $y_i$ to the pre-image of $\Sigma \pitchfork G$ (note that here and below we abuse notation by referring to the pre-images of the tubing arcs in $D^2$ by the same name as their image in $X$). 

Consider a homotopy of the arc $\beta_1$ rel endpoints in the punctured disk $D^2 - \{y_1'\}$ to the arc shown in Figure \ref{fig2}. Since the tube isotopies from Figures \ref{moves2} and \ref{moves} allow $\beta_1$ to be pushed across the pre-image of each double point other than $y_1'$, this gives an isotopy of the tubed surface to one with this new tubing arc (which we continue to label $\beta_1$). This isotopy of $\beta_1$ on the tube diagram, passing $\beta_1$ underneath all other arcs, corresponds to a tube isotopy taking $\mathcal{A}$ to a tubed disk with the tube diagram in Figure \ref{fig2} (recall from Definition \ref{tubeddisk}, crossings between arcs indicate the relative radii of the corresponding tubes). Likewise, the arcs $\beta_2, \dots, \beta_n$ can be isotoped, one at a time and in numerical order, to the purple arcs in Figure \ref{fig3}. Each isotopy on the tube diagram determines a tube isotopy from the tubed disk in Figure \ref{fig2} to the tubed disk in Figure \ref{fig3}, once we require that each $\beta_i$ passes always \emph{over} $\beta_j$ with $j <i$ and \emph{under} $\beta_j$ with $j>i$. 

Now, the tubed disk in Figure \ref{fig3} (isotopic to $D_1$) is in fact the ascending shadow for the inverse of the homotopy shown in Figure \ref{twistdisks} from the disk $D_0$ to the disk $D_1'$. The ascending and descending Whitney arcs for this homotopy are shown in both Figures \ref{fig4} and \ref{twistdisks}. Although from Figure \ref{twistdisks} it is not immediately clear that the disk $D_1'$ is a self-referential disk $(D_0)_\omega$ for some $\omega \in \mathbb{Z}[\pi_1(X)]$, this is shown by the isotopy in Figure \ref{twistiso}. It therefore follows that $D_1$ is isotopic to the self-referential disk $(D_0)_\omega$, with $\mathcal D(\omega) = \Delta(h)$ modulo the image of the homomorphism $\rm{dax}$ in $\mathbb{Z}[\pi_1(X) \setminus 1]$ since $\rm{Dax}(D_0, (D_0)_\omega) = \mathcal D(\omega)$ by Example \ref{Delta}. \qed

\fig{280}{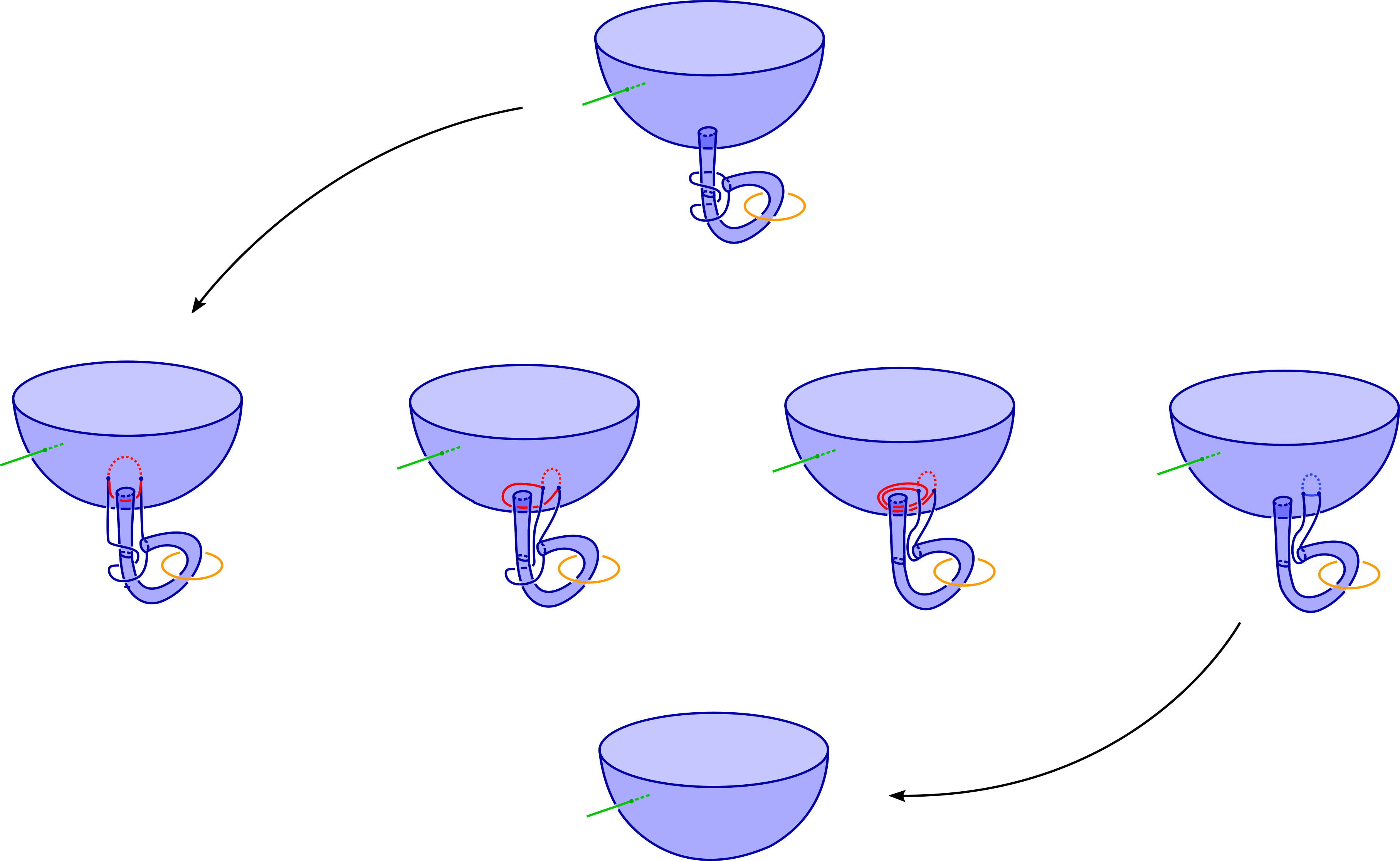}{
\put(-369,224){\small{f.m.}}
\put(-355,110){=}
\put(-230,110){=}
\put(-105,110){=}
\put(-285,10){$D_0$}
\put(-185,210){$D_1'$}
\put(-80,30){\small{w.m.}}
\caption{A regular homotopy from the disk $D_1'$ to $D_0$, consisting of one finger move (f.m.) and one Whitney move (w.m.). The descending Whitney arcs (also featured in Figure \ref{fig4}) are drawn in red, and the ascending Whitney arcs in blue.}
\label{twistdisks}}

\fig{110}{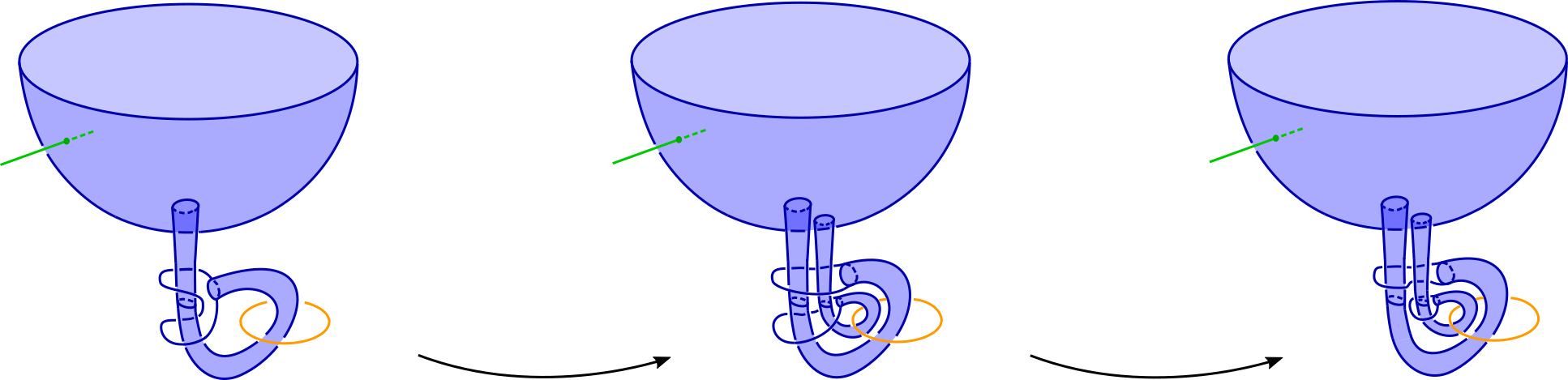}{
\put(-132,-12){\small{\bf Step 2}}
\put(-313,-12){\small{\bf Step 1}}
\put(-417,-10){\small{$D_1'$}}
\put(-13,-10){\small{$(D_0)_\omega$}}
\caption{An isotopy from the disk $D_1'$ to a self-referential disk. {\bf Step 1:} Isotop $D_1'$ by pushing a sheet of the disk along the tube, splitting it into two parallel tubes. {\bf Step 2:} Use Gabai's light bulb lemma \cite[Lemma 2.3]{dave:LBT} to unlink each tube from the end of the other, isotoping the disk obtained in Step 1 to a self-referential disk $(D_0)_\omega$ for some $\omega \in \mathbb{Z}[\pi_1(X)]$.  The light bulb lemma, and hence the isotopy in this figure, relies on the presence of the dual $G$.}
\label{twistiso}}

\remark
We thank an observant referee for pointing out that in the context of Theorem \ref{one}, $\mathrm{dax}(\pi_3(X))$ is necessarily non-trivial. For, if $G$ is $\pi_1$-negligible, then it has an immersed dual $G’$. Therefore, the connected sum $G\#_gG’$ along a tube representing any non-trivial element $g\in\pi_1(X)$ will in general have non-trivial self-intersection invariant $\mu_2$. However, the non-triviality of the self-intersection invariant $\mu_2$ on $\pi_2(X)$ implies that the Dax homomorphism must also be non-trivial, by Theorem B of Kosanović-Teichner \cite{KT}. 

\remark 
Our result was inspired by Theorem $2.5$ of Gabai \cite{dave:LBL}, which states that any pair of homotopic disks with a common dual \emph{in the boundary} of the ambient $4$-manifold can be put into self-referential form with respect to each other. Indeed, Gabai notes in his final remark below the theorem what we observe here: that allowing the dual in the \emph{interior} rather than on the boundary seems to necessitate some $\pi_1$-negligible assumption on the dual sphere. 

\begin{question} \label{q} 
Can some version of the light bulb theorem for disks be stated even when the map induced by inclusion $\pi_1(X-G) \to \pi_1(X)$ is \emph{not} an isomorphism? The proofs of both  Lemma \ref{homotopyrelG} and Theorem \ref{one} given here require this condition; it is unclear if similar results can be proven otherwise. 
\end{question}

\section{Computing the Dax invariant} \label{appendix} 

Our final section is somewhat of an appendix, which can be considered independently. To complement the applications of the Dax invariant from the previous two sections, we present a new geometric method for its computation. 

Let $h: I \times D^2 \to X^4$ be a regular homotopy between embedded disks, where as in Definition \ref{dax}, the disk $D^2$ is parametrized as the product of two unit intervals $J \times K$. There is thus a smooth map $$\widehat h: I \times J \times K \to I \times J \times X,$$ restricting to the identity on the first two components, which (after perturbing $h$ rel boundary if necessary) is an immersion with finitely many transverse double points $p_1, \dots, p_n$ and no triple points. For our alternative geometric computation of the Dax invariant $\Delta(h)$ of the homotopy $h$, we consider a third related map $$h': I \times J \times K \to I \times X$$ given by sending $(t,s,r) \mapsto (t, h_t(s,r))$. This map is also an immersion, but with double curves rather than double points. 

After a second deformation rel boundary of $h$, it can be assumed that all finger moves of the homotopy occur at the same time $t_{min} \in I$ and that all Whitney moves occur simultaneously at some later time $t_{max} \in I$. The pre-image of the double curves of $h'$ is then a link $L \subset I \times J \times K$ equal to the plat closure of a $2\kappa$-stranded braid $\beta$ in as in Figure \ref{boxes}, whose $\kappa$ minima and $\kappa$ maxima with respect to the $I$ coordinate occur at heights $t_{min}$ and $t_{max}$ respectively. It can also be assumed that the projection of the link $L$ to each $I \times J$ plane is a smoothly immersed curve in general position, with some finite number of crossings. 

\begin{definition}
A homotopy $h$ is said to be in \bit{plat position} if its corresponding map $\widehat h: I \times J \times K \to I \times J \times X$ satisfies the requirements of the previous paragraph. 
\end{definition} 

\fig{130}{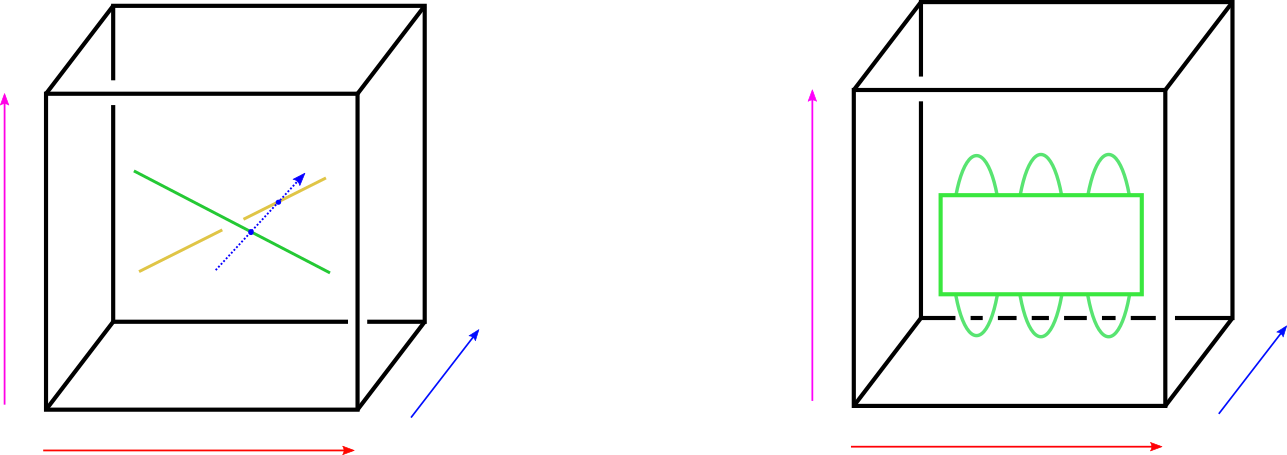}{
\put(-380,55){\small{I}}
\put(-147,55){\small{I}}
\put(-315,-15){\small{J}}
\put(-80,-15){\small{J}}
\put(-75,57){\small{$\beta$}}
\put(-235,16){\small{K}}
\put(-5,16){\small{K}}
\put(-200,65){$\subset$}
\caption{The pre-image of the double curves of the map $h'$ are equal to the plat closure of a $2\kappa$-stranded braid $\beta$ embedded in $I \times J \times K$ (shown on the right). Double points of $\widehat h$ (as shown on the left) occur at each crossing, with respect to the projection in the $K$ direction, between yellow and green strands that are sent to the same arc by $\widehat h$.}
\label{boxes}}

Let $\beta_1, \beta_1', \dots, \beta_\kappa, \beta_\kappa'$ denote the strands of $\beta$, labelled so that $h'(\beta_i)=h'(\beta_i')$ for each $i$. The strands can be oriented so that at each point $x=(t,s,r) \in \beta$ away from the extrema, the tangent space $T_x\beta \subset T_x(I \times J \times K)$ is framed so that the first component of its positive unit vector $\tau_x$ has the same sign as the transverse double point $h_t(s,r) \in X$. Since only oppositely oriented double points are introduced/cancelled during a regular homotopy, this orientation on the strands of $\beta$ extends coherently to one on $L$. 

\begin{definition} \label{crossingdef}
Let $c_{i1}, \dots, c_{i{\ell_i}}$ denote the crossings, with respect to projecting in the $K$ direction onto any $I \times J$ plane, between two strands $\beta_i$ and $\beta_i'$ of the same index. Each crossing $c_{ij}$ can be associated to the homotopy class $g_{ij} \in \pi_1(X)$ of a based loop in $I \times X$ that is the image under $h'$ of oriented paths in $I \times J \times K$ connecting the pre-image $b$ of the basepoint to the double point of the crossing on the understrand (with larger $K$-coordinate), and from the double point of the crossing on the overstrand (with smaller $K$-coordinate) back to $b$. This is the computation for Wall's self-intersection number of an immersion, but done with a ``canonical ordering" of sheets at each double point. Order the crossings so that $g_{ij}$ is non-trivial if and only if $i\leq s$, and let $$\mathcal C(h): = \sum_{i=1}^s \sum_{j=1}
^{\ell_i} \eta_{ij} g_{ij} \in \mathbb{Z}[\pi_1(X) \setminus 1],$$ where the sign $\eta_{ij}=\pm 1$ for $g_{ij} \in \pi_1(X)$ is assigned to using the convention shown in Figure \ref{allcrossings}. 
\end{definition}

\fig{75}{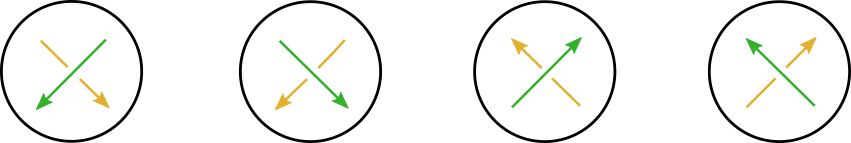}{
\put(-415,-15){$-$}
\put(-290,-15){$+$}
\put(-164,-15){$+$}
\put(-43,-15){$-$}
\caption{The convention for Definition \ref{crossingdef} of assigning a sign $\eta_{ij}=\pm 1$ to each crossing $c_{ij}$ in the projection of the oriented link $L$ in the $K$ direction onto an $I \times J$ plane. Note that these are the only possible crossings, since by construction the strands $\beta_i$ and $\beta_i'$ are oriented in the same direction with respect to the $I$ factor.}
\label{allcrossings}}

\begin{theorem} \label{newway}
If a regular homotopy $h$ is in plat position, then its Dax invariant $\Delta(h)$ is equal to $\mathcal C(h)$ in $\mathbb{Z}[\pi_1(X)\setminus 1]$. 
\end{theorem}

\proof 

Note that the trace $h'$ is equal to the immersion $\pi \circ \widehat h$ for the projection $\pi: I \times J \times X \to I \times X$. It follows that the double points of the immersion $\widehat h$ are in one-to-one correspondence with the crossings of the double curves of $h'$, in the projection to any $I \times J$ plane, between strands $\beta_i$ and $\beta_i'$ of the same index. Compare the left of Figure \ref{boxes} with Figure \ref{daxpic} to see that the corresponding group elements are equal. 

Therefore, all we have left to prove is that for each double point of $\widehat h$ and its corresponding crossing, the associated \emph{signs} given by Definition \ref{dax} and \ref{crossingdef} are equal. To simplify notation below, we drop select subscripts, referring to the double point simply as $p$, and to its pre-images as $x$ and $y$ (points on the understrand and overstrand, respectively).  

Recall, the sign $\epsilon_i$ in Definition \ref{dax} is given by comparing the orientation of the tangent space $T_{p}(I \times J \times X)$ to that of the sum $d \widehat h(T_{x}(I \times J \times K)) \oplus d \widehat h(T_{y}(I \times J \times K))$. To do this, we write a positive ordered bases for each vector space. First, since both the vectors $\tau_x \in T_{x}(I \times J \times K)$ and $\tau_y \in T_{y}(I \times J \times K)$ have non-zero $I$ component, there are vectors $u_x, v_x \in T_{x}(I \times J \times K)$ and $u_y, v_y \in T_{y}(I \times J \times K)$ that are orthogonal to the $I$ direction such that $\{\tau_{x}, u_x, v_x\}$ gives a positive ordered basis for the tangent space $T_{y}(I \times J \times K)$, and $\{\tau_{y}, u_y, v_y\}$ a positive ordered basis for $T_{y}(I \times J \times K)$. 

It follows that $d \widehat h(T_{x}(I \times J \times K)) \oplus d \widehat h(T_{y}(I \times J \times K))$ has positive basis $\{\widehat \tau_{x}, \widehat u_x, \widehat v_x, \widehat \tau_{y}, \widehat u_y, \widehat v_y\}$, where the ``hat" notation is used to denote the image under $d \widehat h$ . Swapping the order of the basis vectors an even number of times obtains a second positive basis $B=\{\widehat \tau_{x}, \widehat \tau_{y}, \widehat u_x, \widehat v_x, \widehat u_y, \widehat v_y\}$. Now, recall that the strands of the braid $\beta$ were oriented so that $\tau_{x}$ has positive $I$ component precisely when $\{ \widehat u_x, \widehat v_x, \widehat u_y, \widehat v_y\}$ is positive ordered basis of $T_{\widehat p}X$, where $\widehat p \subset X$ is the image of $p$ under the projection to $X$. Hence, whether or not the basis $B$ is a positive one for $T_{p}(I \times J \times X) \cong I \times J \times T_{\widehat p}X$ depends only on the sign of the $J$ component of the vector $\tau_y$, which can be determined at each crossing from the orientation and slope of the overstrand, as in Figure \ref{allcrossings}. \qed

\begin{example}\label{exampledax} 
We showcase our new geometric method for computing the Dax invariant by using it to verify Theorem $4.9$ of Gabai \cite{dave:LBL}, discussed earlier in Example \ref{Delta}. For each $g \in \pi_1(X)$, recall that the self-referential homotopy $h_g$ (see Definition \ref{selfref}) from the self-referential disk $D_g$ to the disk $D$ has the Whitney arcs illustrated in both Figures \ref{Dg} and \ref{examplepic}. These arcs ``lift" to the double curves of the map $h'_g$ shown on the right of Figure \ref{examplepic}. By Theorem \ref{newway}, we can diagramatically compute the Dax invariant $\rm{Dax}(D, D_g)$ from these curves, by calculating $\mathcal C(h_g)$ from Definition \ref{crossingdef}. It follows from Figure \ref{examplepic} that $\rm{Dax}(D, D_g)$ is equal to the sum $g + \bar g \in \mathbb{Z}[\pi_1(X) \setminus 1]/\rm{dax}(\pi_3(X))$. 
\end{example}

\fig{280}{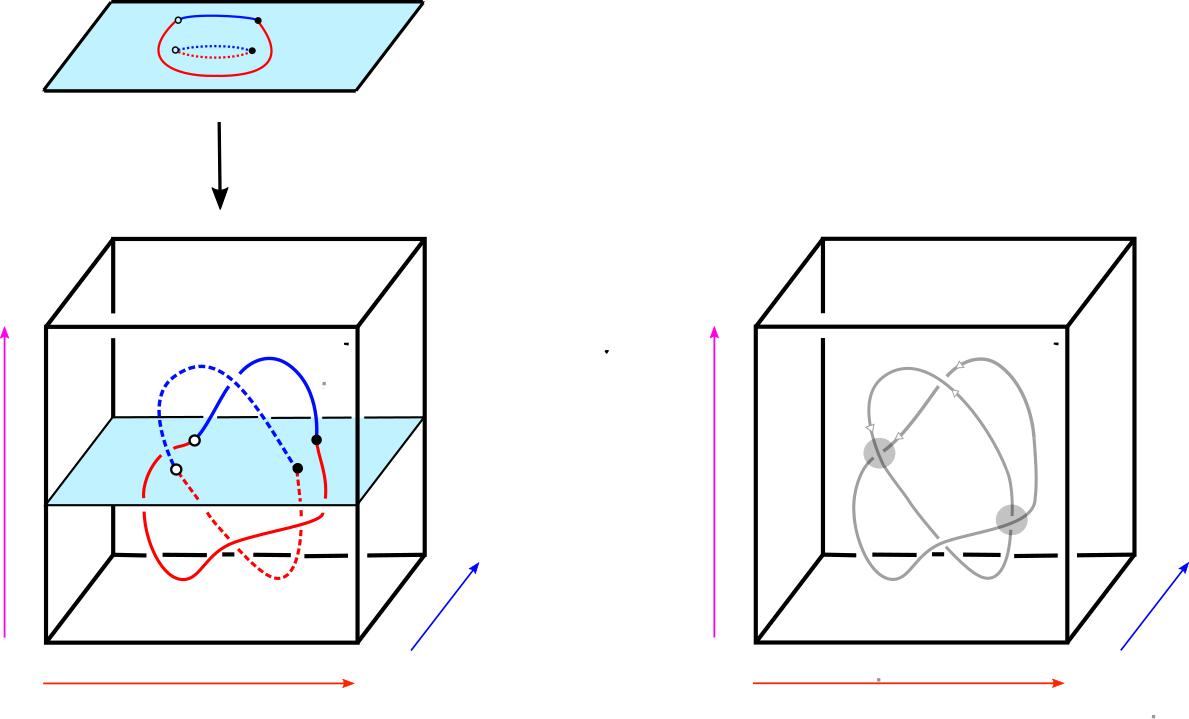}{
\put(-373,215){``lift"}
\put(-280,37){$K$}
\put(-3,37){$K$}
\put(-395,0){$J$}
\put(-120,0){$J$}
\put(-475,85){$I$}
\put(-200,85){$I$}
\put(-250,95){\Large{$=$}}
\caption{The ascending and descending Whitney arcs of the self-referential homotopy $h_g$ (top left), which lift to the double curves of the induced map $h'_g$ (bottom left and right). Note from Figure \ref{Dg} that the homotopy class corresponding to each crossing of the double curves, as in Definition \ref{crossingdef}, is equal to $g$ when the ``outer" curve is the overstrand and $\bar{g}$ when the ``inner" curve is the overstrand. Therefore, the Dax invariant of $h_g$ in this case is equal to sum $g + \bar{g}$ whose terms are determined by comparing each of the highlighted crossings (those with both strands oriented in the same $K$ direction) to those of Figure \ref{allcrossings}.}
\label{examplepic}}

\end{document}